\documentclass[a4paper]{amsart}
\usepackage{latexsym}
\usepackage{mathtools}
\usepackage{amsmath}
\usepackage{amsfonts}
\usepackage{mathrsfs}
\usepackage{slashed}
\usepackage{amscd}
\usepackage{amssymb}
\usepackage{amsfonts}
\usepackage{float}
\usepackage{esint}
\usepackage{caption}
\usepackage{hyperref}
\usepackage{tikz-cd}

\usepackage[utf8]{inputenc}
\usepackage[T1]{fontenc}   
\usepackage{graphicx}
\usepackage{color}
\allowdisplaybreaks 
\definecolor{Blue}{rgb}{0.,0.,1.}
\definecolor{Red}{rgb}{1.,0.,0.}
\definecolor{Green}{rgb}{0.,1.,0.}

\let\origmaketitle\maketitle
\def\maketitle{
  \begingroup
	\origmaketitle
  \endgroup
	}	
	
\makeatletter
\renewcommand{\author}[2][]{%
  \def\@tempa{#1}
  \ifx\@empty\authors
    \ifx\@tempa\@empty
      \gdef\shortauthors{#2}%
    \else
      \gdef\shortauthors{#1}%
    \fi
    \gdef\authors{\author{#2}}%
  \else
    \ifx\@tempa\@empty
      \g@addto@macro\shortauthors{\and#2}%
    \else
      \g@addto@macro\shortauthors{\and#1}%
    \fi
    \g@addto@macro\authors{\and\author{#2}}%
  \fi
}
\renewcommand{\address}[2][]{\g@addto@macro\authors{\address{#1}{#2}}}
\def\@setauthors{%
  \begin{center}%
    \footnotesize
    \vspace{20pt}
    \let\and\@empty
    \def\author##1{\advance\@tempcnta\@ne}%
    \def\address##1##2{\advance\@tempcntb\@ne}%
    \@tempcnta=\z@  \@tempcntb=\z@
    \authors
    \ifnum\@tempcnta>\@ne \ifnum\@tempcntb=\@ne
        \oneaddress
      \else
        \sepaddresses
      \fi
    \else
      \oneaddress
    \fi
  \end{center}%
}
\def\oneaddress{%
  \begingroup
  \let\author\@iden \let\address\@gobbletwo
  \renewcommand{\andify}{%
    \nxandlist{\unskip, }{\unskip{} and~}{\unskip, and~}}%
  \uppercasenonmath\authors
  \andify\authors
  \authors
  \endgroup
  \begingroup \let\and\relax \let\author\@gobble
  \def\address##1##2{\unskip\\[10pt] \itshape##2}%
  \authors
  \endgroup
}
\def\sepaddresses{%
  \begingroup
    \baselineskip10\p@\relax
    \def\address##1##2{ ({\itshape##2}\/)}
    \def\author##1{\def\temp{##1}\leavevmode\uppercasenonmath\temp\temp}%
    \nxandlist
      {,\\[\baselineskip]}
      {\\[\baselineskip] \textsc{\lowercase{and}}\\[\baselineskip]}
      {,\\[\baselineskip]\textsc{\lowercase{and}}\\[\baselineskip]}
      \authors 
    \authors
  \endgroup
}
\def\maketitle{\par
  \@topnum\z@
  \@setcopyright
  \thispagestyle{firstpage}%
  \uppercasenonmath\shorttitle
  \ifx\@empty\shortauthors \let\shortauthors\shorttitle
  \else
    \newcommand{\@xuppercasenonmath}[1]{\toks@\@emptytoks
      \@xp\@skipmath\@xp\@empty##1$$%
      \edef##1{\@nx\protect\@nx\@upprep\the\toks@}}%
    \@xuppercasenonmath\shortauthors
    \def\@@and{AND}
    \renewcommand{\andify}{%
      \nxandlist{\unskip, }{\unskip{ }\@@and{ }}{\unskip, \@@and{ }}}%
    \andify\shortauthors
  \fi
  \@maketitle@hook
  \begingroup
  \@maketitle
  \endgroup
  \c@footnote\z@
  \@cleartopmattertags
}
\def\@maketitle{%
  \normalfont\normalsize
  \let\@makefntext\noindent
  \@adminfootnotes
  \ifx\@empty\addresses\else \@footnotetext{\@setotheraddresses}\fi
  \global\topskip68\p@\relax
  \@settitle
  \ifx\@empty\authors \else \@setauthors \fi
  \ifx\@empty\@dedicatory
  \else
    \baselineskip26\p@
    \vtop{\centering{\footnotesize\itshape\@dedicatory\@@par}%
      \global\dimen@i\prevdepth}\prevdepth\dimen@i
  \fi
  \toks@\@xp{\shortauthors}\@temptokena\@xp{\shorttitle}%
  \edef\@tempa{\@nx\markboth{\the\toks@}{\the\@temptokena}}\@tempa
  \@setabstract
  \normalsize
  \if@titlepage
    \newpage
  \else
    \dimen@34\p@ \advance\dimen@-\baselineskip
    \vskip\dimen@\relax
  \fi
} 
\renewcommand{\thanks}[1]{%
  \ifx\@empty\thankses
    \gdef\thankses{\thanks{#1}}%
  \else
    \g@addto@macro\thankses{\endgraf\thanks{#1}}%
  \fi}
\def\@setthanks{\def\thanks##1{\noindent##1\@addpunct.}\thankses}
\renewcommand{\curraddr}[2][]{%
  \ifx\@empty\addresses
    \gdef\addresses{\curraddr{#1}{#2}}%
  \else
    \g@addto@macro\addresses{\endgraf\curraddr{#1}{#2}}%
  \fi}
\renewcommand{\email}[2][]{%
  \ifx\@empty\addresses
    \gdef\addresses{\email{#1}{#2}}%
  \else
    \g@addto@macro\addresses{\endgraf\email{#1}{#2}}%
  \fi}
\renewcommand{\urladdr}[2][]{%
  \ifx\@empty\addresses
    \gdef\addresses{\urladdr{#1}{#2}}%
  \else
    \g@addto@macro\addresses{\endgraf\urladdr{#1}{#2}}%
  \fi}
\def\@setotheraddresses{%
  \def\curraddr##1##2{\noindent
    \emph{Current address\@ifnotempty{##1}{ of ##1}}:\space
      ##2\@addpunct.}%
  \def\email##1##2{\noindent
    \emph{E-mail address\@ifnotempty{##1}{ of ##1}}:\space
      \texttt{##2}}%
  \def\urladdr##1##2{\noindent
    \emph{WWW address\@ifnotempty{##1}{ of ##1}}:\space
      \texttt{##2}}%
  \addresses
}
\let\enddoc@text\relax
\makeatother

\newcounter{smallarabics}
\newenvironment{arabicenumerate}
{\begin{list}{{\normalfont\textrm{(\arabic{smallarabics})}}}
  {\usecounter{smallarabics}\setlength{\itemindent}{0cm}
   \setlength{\leftmargin}{5ex}\setlength{\labelwidth}{4ex}
   \setlength{\topsep}{0.75\parsep}\setlength{\partopsep}{0ex}
   \setlength{\itemsep}{0ex}}}
{\end{list}}

\newcounter{smallroman}

\newcommand{\ben}{\begin{arabicenumerate}}  
\newcommand{\een}{\end{arabicenumerate}}

\def\init{\setcounter{equation}{0}}

\newtheorem{theoreme}{Theorem }[section]
\newtheorem{proposition}[theoreme]{Proposition}

\newtheorem{lemma}[theoreme]{Lemma}
\newtheorem{definition}[theoreme]{Definition}

\newtheorem{remark}[theoreme]{Remark}
\newtheorem{example}[theoreme]{Example}
\newcommand{\beq}{\begin{equation}}
\newcommand{\eeq}{\end{equation}}

\newcommand{\bex}{\begin{example}}
\newcommand{\eex}{\end{example}}
\def\bel{\begin{lemma}}
\def\eel{\end{lemma}}
\def\bet{\begin{theoreme}}
\def\eet{\end{theoreme}}
\def\bed{\begin{definition}}
\def\eed{\end{definition}}
\def\ber{\begin{remark}}
\def\eer{\end{remark}}

\def\rr{{\mathbb R}}

\def\cc{{\mathbb C}}
\def\nn{{\mathbb N}}

\def\bar{\overline}

\def\cinf{C^\infty}
\def\proof{
\noindent{\bf Proof.}\ \ }

\def\cY{{\mathcal Y}}

\def\cL{{\mathcal L}}
\def\cS{{\mathcal S}}

\def\cD{{\mathcal D}}

\def\cW{{\mathcal W}}

\def\i{{\rm i}}
\def\qed{$\Box$\medskip}
\def \p{ \partial}
\def\12{\frac{1}{2}}
\def\14{\frac{1}{4}}

\def\bbbone{{\mathchoice {\rm 1\mskip-4mu l} {\rm 1\mskip-4mu l}
{\rm 1\mskip-4.5mu l} {\rm 1\mskip-5mu l}}}
\def\one{\bbbone}
\def\cH{{\mathcal H}}

\def\coinf{C_0^\infty}

\def\cT{{\mathcal T}}
\def\cF{{\mathcal F}}

\def\cX{{\mathcal X}}

\def\S{{\mathcal S}}

\def \p{ \partial}
\def\12{\frac{1}{2}}
\def\e{{\rm e}}

\def\Op{{\rm Op}}

\newcommand{\traa}[1]{\mskip-6mu\upharpoonright_{#1}}

\def\cE{{\mathcal E}}

\def\WF{{\rm WF}}
\makeatletter
\newcommand*{\defeq}{\mathrel{\rlap{%
                     \raisebox{0.3ex}{$\m@th\cdot$}}%
                     \raisebox{-0.3ex}{$\m@th\cdot$}}%
                     =}
\makeatother
\makeatletter
\newcommand*{\eqdef}{=\mathrel{\rlap{%
                     \raisebox{0.3ex}{$\m@th\cdot$}}%
                     \raisebox{-0.3ex}{$\m@th\cdot$}}%
                     }
\makeatother

\def\Sol{{\rm Sol}_{\rm sc}}

\def\cinfb{C^{\infty}_{\rm b}}
\def\rx{{\rm x}}
\def\ry{{\rm y}}

\DeclareMathOperator{\Ker}{Ker}
\DeclareMathOperator{\Dom}{Dom}
\DeclareMathOperator{\supp}{supp}

\def\dual{\!\cdot \!}

\def\cB{\mathcal{B}}

\def\zero{{\mskip-4mu{\rm\textit{o}}}}

\def\cN{{\mathcal N}}

\def\BT{{ BT}}

\def\tg{{\tilde g}}
\def\tnab{{\tilde \nabla}}

\def\oibis{{\rm out/in}}
\def\oito{\pm\infty}
\def\outin{\pm\infty}
\def\oi{{\rm out/in}}

\def\tosim{\xrightarrow{\sim}}

\def\zero{{\mskip-4mu{\rm\textit{o}}}}
\def\Cl{{\rm Cliff}}

\newcommand{\bea}{\begin{aligned}}
\newcommand{\beal}{\begin{array}{l}}
\newcommand{\eeal}{\end{array}}
\newcommand{\eea}{\end{aligned}}

\def\SO{{\rm SO}}
\def\Spin{{\rm Spin}}

\def\th{\tilde{h}}
\def\CAR{{\rm CAR}}
\def\spinup{\Spin^{\uparrow}}\def\soup{\SO^{\uparrow}}
\def\blibli{\tilde}
\def\maketitle{
  \begingroup
  \def\uppercasenonmath##1{} 
  \let\MakeUppercase\relax 
	\origmaketitle
  \endgroup
	}
\begin{document}
\pagestyle{plain}

\title{\large Hadamard property of the \emph{in} and \emph{out} states \\ for Dirac fields on asymptotically static spacetimes
}
\author{\normalsize Christian \textsc{G\'erard}}
\address{Universit\'e Paris-Saclay, D\'epartement de Math\'ematiques, 91405 Orsay Cedex, France}
\email{christian.gerard@math.u-psud.fr}
\date{March 2021}
\author{\normalsize Th\'eo \textsc{Stoskopf} }
\email{theo.stoskopf@universite-paris-saclay.fr}
\begin{abstract} We consider massive Dirac equations on asymptotically static spacetimes with a Cauchy surface of bounded geometry.  We prove that the associated quantized Dirac field admits {\em in } and {\em out} states, which are asymptotic vacuum states when some time coordinate tends to $\mp\infty$. We also show that the {\em in /out} states are Hadamard states. 
\end{abstract}

\maketitle

\section{Introduction}\init\label{sec-1}
\subsection{{\em In/out} vacuum states}
The construction of a {\em distinguished quantum state} for a quantized field on a curved background   has been studied extensively in various contexts in Quantum Field Theory.   

If the background spacetime has no global symmetries but only asymptotic ones,  one can try to specify a distinguished quantum state by its asymptotic behavior, for example at early or late times.

An often studied situation arises when the background spacetime $(M, g)$ has a product structure $M = \rr\times \Sigma$ and the metric $g$ becomes asymptotic to {\em static} metrics when $t\to \pm \infty$. One can then at least heuristically consider {\em asymptotic vacua},  the so-called \textit{in} and \textit{out} states, which look like vacuum states for the asymptotic static metrics when $t\to \mp\infty$.

Let us mention for example the  wave and Klein-Gordon fields  on Minkowski space, in external electromagnetic potentials \cite{isozaki,rui,seiler}, or on curved spacetimes with special a\-symp\-totic symmetries, \cite{wald0,DK0,DK1,DK2}. 

Besides the existence of the \textit{in} and \textit{out} states, an important question is to ensure that  they 
 satisfy the \emph{Hadamard condition} \cite{KW}. 
 
 Nowadays regarded as an indispensable ingredient in the perturbative construction of interacting fields (see e.g.  the recent reviews \cite{HW,FV2}), this property accounts for the correct short-distance behaviour of two-point functions. It can be conveniently formulated as a condition on the \emph{wave front set} of the state's two-point functions \cite{radzikowski}.

The above questions were solved  in \cite{GW} for {\em massive Klein-Gordon fields}, using a combination of scattering theory arguments and global pseudodifferential calculus. 
 
 In this paper we  consider this problem for {\em massive Dirac fields}, using similar methods. Let us now describe  in more details the results of this paper.
 \subsection{Results}
 \subsubsection{Asymptotically static spacetimes}
We  will consider a  spacetime of {\em even dimension}  $n$ of the form $M= \rr\times \Sigma$, where $\Sigma$ is a $d$-dimensional manifold, equipped with a metric
 \[
 g= - c^{2}(x)dt^{2}+ (d\rx^{i}+ b^{i}(x)dt)h_{ij}(x)(d\rx^{j}+ b^{j}(x)dt),
 \]
 where $x=(t, \rx)\in M$, $c\in \cinf(M; \rr)$ is a strictly positive function, $b\in \cinf(M; T\Sigma)$ and $h\in \cinf(M; \otimes^{2}_{\rm s}T^{*}\Sigma)$ is a $t$-dependent Riemannian metric on $\Sigma$. We will assume that when $t\to \pm \infty$ the metric $g$ converges to {\em static metrics}
 \[
 g_{\oi}=- c_{\oi}(\rx)dt^{2}+ h_{\oi}(\rx)d\rx^{2}.
 \]
 The convergence of $g$ to $g_{\oi}$ is assumed to be uniform in the space variable $\rx$.  More precisely, one assumes that there exists $\mu>0$ such that
 \beq\label{convero}
\left. \begin{array}{l}
  \p_{t}^{k}\p_{\rx}^{\alpha}(h(x)- h_{\oi}(\rx))\in O(\langle t\rangle^{-\mu-k }), \\[2mm]
    \p_{t}^{k}\p_{\rx}^{\alpha}b(x)\in O(\langle t\rangle^{-1-\mu-k }), \\[2mm]
    \p_{t}^{k}\p_{\rx}^{\alpha}(c(x)- c_{\oi}(\rx))\in O(\langle t\rangle^{-\mu-k }),
 \end{array}\right. \ k\in \nn, \alpha\in \nn^{d}, 
 \eeq
 in an appropriate uniform sense in $\rx\in \Sigma$, using the notion of Riemannian manifolds of {\em bounded geometry}, see  hypotheses ${\rm (H1)}$, ${\rm (H2)}$ in \ref{sec1.1.1}  for the precise  formulation.

 Note that  the existence of the in/out state is usually deduced from a result of existence and completeness of {\em M\"{o}ller operators} for classical fields which typically requires a short-range condition   $\mu>1$ in \eqref{convero}. As  was the case in \cite{GW},  only the weaker  condition $\mu>0$  is needed in our paper.
\subsubsection{Dirac operators}
 We consider a Dirac operator 
  \[
  D= \slashed{D}+ m
  \]
   and assume 
 that $m(t, \rx)$ converges to $m_{\oi}(\rx)$ when $t\to \pm\infty$ in a similar uniform way, see hypothesis ${\rm (H3)}$ in \ref{sec1.1.1}  for the precise  formulation.

  It follows  that   $D$  converges when $t\to \pm\infty$ to {\em asymptotic Dirac operators}  $D_{\oi}$, which are associated to the static metrics $g_{\oi}$.

The vector field $\p_{t}$ is  Killing for the static metrics $g_{\oi}$, which implies that one can  define the {\em vacuum states} $\omega^{\rm vac}_{\oi}$ for $D_{\oi}$, see Subsect. \ref{secvac}, using the projections
\[
c^{\pm{\rm vac}}_{\oi}\defeq \one_{\rr^{\pm}}(H_{\oi}), 
\]
where $H_{\oi}$ are selfadjoint operators on $L^{2}(\Sigma; S(\Sigma))$, for the  canonical Hilbertian scalar product on $S(\Sigma)$.  The operators  $H_{\oi}$ are the {\em generators} of the unitary group induced by  the spinorial Lie derivative $\cL_{\p_{t}}$ on solutions of $D_{\oi}\psi=0$, see Subsect. \ref{secvac}.

To define the vacuum states $\omega^{\rm vac}_{\oibis}$ in an unambiguous way, one needs to assume that    
   \beq\label{kero}
\Ker H_{\oi}= \{0\},
  \eeq 
ie the absence of {\em zero modes}. If \eqref{kero} is violated, then in physics language one needs to decide  if zero modes are considered as particles or as anti-particles. 
  
In this paper, we strengthen  \eqref{kero} by requiring that 
\[
0\not\in \sigma(H_{\oi}),
\]
 see hypothesis $({\rm H}4)$, ie that the asymptotic Dirac operators $D_{\oi}$ are {\em massive} in the terminology of \ref{energy-gap}.

\subsubsection{Existence of the {\em in/out} states}
Let us now explain the definition of the {\em in/out} states for $D$. We set $\Sigma_{s}= \{s\}\times \Sigma$ and fix the reference time $t=0$.

Denoting by $U(t,s): \coinf(\Sigma_{s}; S(\Sigma_{s}))\to \coinf(\Sigma_{t}; S(\Sigma_{t}))$ the {\em Cauchy evolution} for the Dirac operator $D$ one  expects that the limits
\begin{equation}
\label{limitoto}
c^{\pm}_{\oi}= \lim_{t\to \pm\infty}U(0, t)c^{\pm{\rm vac}}_{\oi}U(t, 0)
\end{equation}
exist in an appropriate sense.   
The  $c^{\pm}_{\oi}$ are supplementary projections acting on the space of Cauchy data at time $t=0$, which are  selfadjoint for the canonical Hilbertian scalar product. Therefore one can associate to  $c^{\pm}_{\oi}$ quasi-free states for the free Dirac field on $M$. Concretely one considers  pairs of operators
\[
\Lambda^{\pm}_{\oibis}: \coinf(M; S(M))\to \cinf(M; S(M))
\]
defined by
\beq\label{turlu}
\Lambda^{\pm}_{\oibis}(t, s)= U(t, 0)\i \gamma(n)c^{\pm}_{\oi}U(0, s),
\eeq
where we write $\Lambda^{\pm}_{\oi}$ as operator-valued Schwartz kernels in the time variable,  ie we use the formal identity
\[
Au(t)= \int_{\rr}A(t,s)u(s)ds,
\]
to define the 'time kernel' of some operator $A$ acting on $M$. In \eqref{turlu}
$n$ is the future directed unit normal to $\Sigma_{0}$ and  $\gamma$ are the 'gamma matrices' (or Clifford multiplications)  obtained from the spin structure on  $(M, g)$.

The operators  $\Lambda^{\pm}_{\oi}$
 satisfy:
  \[
  \begin{array}{rl}
  {\rm (i)}& \Lambda_{\oi}^{\pm}\geq 0,\\[2mm]
{\rm (ii)}&\Lambda_{\oi}^{+}+ \Lambda_{\oi}^{-}= \i G,\\[2mm]
{\rm (iii)}&D\circ \Lambda_{\oi}^{\pm}= \Lambda_{\oi}^{\pm}\circ D=0,
\end{array}
\]
where $G= G_{\rm ret}- G_{\rm adv}$ is the causal propagator for $D$ and the positivity in ${\rm (i)}$ is formulated with respect to the canonical indefinite Hermitian form on $\coinf(M; S(M))$ for which $D$ is formally selfadjoint.  It follows that   $\Lambda^{\pm}_{\oi}$ are {\em two-point functions}  of  the sought after {\em in/out} states $\omega_{\oibis}$ for $D$. 

\subsubsection{Hadamard property of the {\em in /out} states}
As explained above, the Hadamard condition allows to select among the plethora of states  the physically meaningful ones, which should resemble the Minkowski vacuum, at least in the vicinity of any point of $M$.

The microlocal definition of Hadamard states for Dirac fields was first introduced by Hollands \cite{Ho}, who also proved its equivalence with the older characterization by the short distance asymptotics of its two-point functions. Hadamard states for Dirac fields were further studied in  \cite{K, Kr, SV, S}.  

To our knowledge the first paper proving existence of Hadamard states for Dirac fields in the general case is the recent paper by Islam and Strohmaier \cite{IS}, although  the construction of Hadamard states  by  the deformation argument of Fulling, Narcowich and Wald \cite{FNW} was quite probably known to experts.  

Another construction of Hadamard states on spacetimes of bounded geometry was given in \cite{St} using global pseudodifferential calculus on a Cauchy surface. The methods used in the present paper are to a large extend an adaptation of the strategy in \cite{St} to a scattering situation.

Let us now state the main result of this work, referring the reader to  Subsect. \ref{sec1.1} for hypotheses ${\rm (H)}$.
\begin{theoreme}\label{mainmain}
Assume hypotheses ${\rm (Hi)}$, $1\leq i\leq 4$. Then
\ben
\item The  norm limits \eqref{limitoto} exist and define by \eqref{turlu} {\em pure} quasi-free states $\omega_{\oibis}$ called  the {\em in/out  vacuum states}.
\item  $\omega_{\oibis}$ is a {\em Hadamard state}, ie 
\[
\WF(\Lambda^{\pm}_{\oibis})\subset \cN^{\pm}\times \cN^{\pm}
\]
where $\cN^{\pm}$ are the two connected components of the characteristic set $\cN= \{(x, \xi)\in T^{*}M\setminus \zero: \xi\dual g^{-1}(x)\xi= 0\}$ of $D$. 
\een
\end{theoreme}
\subsection{Outline of the proof} 
Let us now briefly explain the main ingredients in the proof of Thm. \ref{mainmain}, which follows the general strategy in \cite{St}. The first step consists in reducing the metric to the simpler form 
\[
g= - dt^{2}+ h(t, \rx)d\rx^{2},
\]
where the  time-dependent Riemannian metric  $h(t, \rx)d\rx^{2}$ on $\Sigma$  converges to  Riemannian metrics $h_{\oi}(\rx)d\rx^{2}$ when $t\to \pm\infty$. This is done in the usual way, by combining a conformal transformation and the well-known argument using the flow of the vector field $\nabla t$.

One can use the  covariance of  Dirac operators  and two-point functions under conformal transformations, see  Subsect. \ref{sec0.2} and \ref{sec5.1b.1}, to reduce ourselves to this simple situation.

In a second step, one uses parallel transport with respect to  the vector field $ \p_{t}$ to identify the spinor bundles at different times, and to reduce the Dirac equation
$D\psi=0$ to a time-dependent Schroedinger equation:
\[
\p_{t}\psi- \i H(t)\psi=0,
\]
where  $H(t)= H(t, \rx, \p_{\rx})$ is a first order elliptic differential operator on $\Sigma$.

The third step is analogous to \cite{St}, where Hadamard states for Dirac fields are constructed using pseudodifferential calculus, with the difference that in our case we need to control the behavior of various operators when $t\to \pm\infty$.

 We construct time-dependent projections $P^{\pm}(t)$ such that 
 \[
 \begin{array}{rl}
(1)& P^{\pm}(t)- \one_{\rr^{\pm}}(H_{\oi})\in O(t^{-\mu})\hbox{ when }t\to \pm\infty,\\[2mm]
 (2)&\p_{t}U(0, t)P^{\pm}(t)U(t, 0)\in O(t^{-1-\mu})\Psi^{-\infty},
\end{array}
  \]
 where $\Psi^{-\infty}$ is some ideal of smoothing operators on $\Sigma$.   (1) implies that to prove  the existence of the limits \eqref{limitoto}, it sufices to consider instead
 \[
 \lim_{t\to \pm\infty}U(0, t)P^{\pm}(t)U(t, 0)
 \]
 which exists by (2) and the Cook argument. This prove the existence of  the {\em in /out} states $\omega_{\oibis}$ for $D$. Integrating (2) from $0$ to $\pm\infty$, we also obtain that 
 \[
 c^{\pm}_{\oi}- P^{\pm}(0)\hbox{ are smoothing operators on }\Sigma.
 \]
It is shown in  \cite{St} that $P^{\pm}(0)$ are projections which generate a Hadamard state, which, since $c^{\pm}_{\oi}- P^{\pm}(0)$ are smoothing, proves the Hadamard property of $\omega_{\oibis}$.

 \subsection{Plan of the paper}
 Let us now discuss the plan of this paper.
 
 In Sect. \ref{sec0} we recall the quantization of Dirac fields on curved spacetimes. 
 In Sect. \ref{sec1} we describe the geometric framework of asymptotically static spacetimes and the spin structures and Dirac operators on such spacetimes.
 
 In Sect. \ref{sec2} we give a brief overview of Shubin's   pseudodifferential calculus on manifolds of bounded geometry and of its time-dependent version that we will use in this paper. 
 Finally Sect. \ref{sec4} contains the proof of Thm. \ref{mainmain} and the various reduction procedures that are used.
\subsection{Notation}

\subsubsection{Lorentzian manifolds}
We use the mostly $+$  signature convention for Lorentzian metrics. All Lorentzian manifolds considered in this paper will be {\em orientable} and connected.

\subsubsection{Bundles}
If  $E\xrightarrow{\pi}M$ is a bundle we denote by $\cinf(M; E)$ resp. $\coinf(M; E)$ the set of smooth resp. smooth and compactly supported sections of $E$.

If $E\xrightarrow{\pi}M$ is a vector bundle we denote by $\cD'(M; E)$ resp. $\cE'(M; E)$ the space of distributional resp. compactly supported distributional sections of $E$.
\subsubsection{Matrices}
Since we will often use frames of vector bundles we will  denote by $\pmb{M}$ a matrix in  $M_{n}(\rr)$ or $M_{N}(\cc)$  and by $M$ the associated endomorphism. 
 \subsubsection{Frames and frame indices} We use the letters $0\leq a \leq d$ for frame indices on $TM$ or $T^{*}M$, and $1\leq a\leq d$ for frame indices on $T\Sigma$ or $T^{*}\Sigma$, if $\Sigma\subset M$ is a space like hypersurface.  If $g$ is a metric on $M$ and $(e_{a})_{0\leq a\leq d}$ is a local frame of $TM$ we set $g_{ab}= e_{a}\dual g e_{b}$ and $g^{ab}= e^{a}\dual g^{-1}e^{b}$, where $(e^{a})_{0\leq a \leq d}$ is the dual frame. 
 
  We use capital letters $1\leq A\leq N$ for frame indices of the spinor bundle $S(M)$.
  
  If $\cF$ is for example a local frame of $TM$  we denote by $\cF\pmb{t}$ the frame obtained by the right action of $\pmb{t}\in M_{n}(\rr)$ on $\cF$.

  We use capital letters $1\leq A\leq N$ for frame indices of the spinor bundle $S(M)$.

\subsubsection{Vector spaces}
if $\cX$ is a real or complex vector space, we denote by $\cX'$ its dual. If  $\cX$ is a complex vector space we denote by $\cX^{*}$ its anti-dual, ie the space of anti-linear forms on $\cX$ and by  $\bar{\cX}$ its conjugate, ie $\cX$ equipped with the complex structure $-\i$. 

A linear map $a\in L(\cX, \cX')$ is a bilinear  form on $\cX$, whose action on pairs of vectors is denoted  by $x_{1}\dual a x_{2}$. Similarly a linear map $a\in L(\cX, \cX^{*})$ is a sesquilinear form on $\cX$, whose action is denoted by $\bar{x}_{1}\dual a x_{2}$. We denote by $a'$, resp. $a^{*}$ the transposed resp. adjoint of $a$.  The space of symmetric resp. hermitian forms on $\cX$ is denoted by $L_{\rm s}(\cX, \cX')$ resp. $L_{\rm h}(\cX, \cX^{*})$.
\subsubsection{Maps}
We write $f:A\xrightarrow{\sim}B$ if $f: A\to B$ is a bijection. We use the same notation if $A, B$ are topological spaces resp. smooth manifolds, replacing bijection by homeomorphism, resp. diffeomorphism.

 \section{Quantization of Dirac equations on curved spacetimes}\label{sec0}\init
 In this section we recall well-known facts, see eg \cite{Di, LM, Ho, T} about Dirac equations and Dirac quantum fields on curved spacetimes.

\subsection{Dirac equations on curved spacetimes}
Let us denote by $\soup(1, d)$ and $\spinup(1, d)$ the restricted Lorentz and Spin groups (ie the connected component of $Id$ in ${\rm O}(1, d)$ and ${\rm Pin} (1, d)$) and $Ad: \spinup(1, d)\to \soup(1, d)$ the double sheeted covering.

We recall that a {\em spacetime} is  an oriented and time oriented Lorentzian manifold.

\subsubsection{Spin structures}
Let $(M, g)$ a  spacetime of even dimension $n= 1+d$ and let $P\soup(M, g)$ the $\soup(1, d)$-principal bundle of oriented and time oriented orthonormal frames of $TM$.

  We recall that a {\em spin structure} on $(M, g)$ is given by a $\spinup(1,d)$-principal bundle $P{\rm Spin}(M, g)$ with a bundle morphism  $\chi: P{\rm Spin}(M, g)\to P\soup(M, g)$ such that  the following diagram commutes:
  \beq\label{e1.-3}
 \begin{tikzcd}
\Spin^{\uparrow}(1,d) \arrow[r] \arrow[dd, "Ad"] &P\Spin(M, g) \arrow[rd, "\pi'"] \arrow[dd, "\chi"] &   \\
                   &                         & M. \\
\SO^{\uparrow}(1, d) \arrow[r]            & P\SO^{\uparrow}(M, g) \arrow[ru, "\pi"]            &  
\end{tikzcd}
\eeq

We assume that $(M, g)$ has a spin structure $P\Spin(M, g)$. Let us recall that a 
 Lorentzian manifold admits  a spin structure if and only if  its  second Stiefel--Whitney class $w_{2}(TM)$ is trivial, see \cite{Mi, Na}. It admits a {\em unique} spin structure if in addition  its  first Stiefel--Whitney class $w_{1}(M)$ is trivial, which is equivalent to the fact that $M$ is orientable, see e.g.~\cite{Na}. 
  In our situation,   $M$ is orientable hence spin structures on $(M, g)$ are unique  if they exist.  If $n=4$ and $(M, g)$ is globally hyperbolic it admits a (unique) spin structure, see \cite{spinors1, spinors2}.
  
We  denote by $\Cl(M, g)$, $S(M)$ the associated Clifford and spinor bundles. 

The map $TM\to End(S(M))$ obtained from the embedding $TM\to \Cl(M, g)$ and the canonical map $\Cl(M, g)\to End(S(M))$ will be denoted by
\beq\label{troup.-1}
TM\ni u\mapsto \gamma(u)\in End(S(M)),
\eeq
 and is often called the {\em Clifford multiplication}. The spin connection will be denoted by $\nabla^{S}$.

It is well known see eg \cite{T}, \cite[Sect. 17.6]{G} that there exists a (unique up to multiplication by strictly positive constants) non degenerate Hermitian form $\beta$ acting on the fibers of $S(M)$ such that
\beq\label{troup.0}
\begin{array}{l}
\gamma^{*}(u)\beta = - \beta \gamma(u), \ u\in TM, \\[2mm]
 \i \beta \gamma(e)>0, \hbox{for all  }e\in TM \hbox{ time-like and future directed},\\[2mm]
 u\dual \bar{\psi}\dual \beta \psi= \overline{\nabla^{S}_{u}\psi}\dual\beta \psi+ \bar{\psi}\dual \beta \nabla^{S}_{u}\psi, \ \forall
  u\in \cinf(M; TM), \psi\in \cinf(M; S(M)).
\end{array}
\eeq
For later use we summarize the properties of $\nabla^{S}$, $\gamma$ and $\beta$ that we will need. We have:
\beq\label{troup.1}
\begin{array}{l}
\nabla_{u}^{S} \gamma(v)\psi= \gamma(v)\nabla_{u}\psi+ \gamma(\nabla_{u}v)\psi, \\[2mm]
u\dual \bar{\psi}\dual \beta \psi= \overline{\nabla^{S}_{u}\psi}\dual\beta \psi+ \bar{\psi}\dual \beta \nabla^{S}_{u}\psi,\\[2mm]
 u, v\in \cinf(M; TM), \psi\in \cinf(M; S(M)),
\end{array}
\eeq
where $\nabla$ is the metric connection on $(M, g)$
\subsubsection{Dirac operators}

Fixing a smooth section $m\in \cinf(M; L(S(M)))$ with $m^{*}\beta= \beta m$, we consider a Dirac operator 
\beq\label{defdedir}
D= \slashed{D}+ m,
\eeq
where $\slashed{D}$ is locally expressed (on  an open set $U\subset M$ over which $S(M)$ and $TM$ are trivialized) as
\[
\slashed{D}= g^{ab}\gamma(e_{a})\nabla^{S}_{e_{b}},
\]
where $(e_{a})_{0\leq a\leq d}$ is a local frame over $U$.

\subsubsection{Selfadjointness}
For $\psi_{1}, \psi_{2}\in \cinf(M; \S(M))$ one defines   the $1$-form $J(\psi_{1}, \psi_{2})$ $\in \cinf(M; T^{*}M)$ by 
\[
J(\psi_{1}, \psi_{2})\dual u\defeq \overline{\psi}_{1}\dual \beta \gamma(u)\psi_{2}, \quad u\in \cinf(M; TM),
\]
and one deduces from \eqref{troup.1} that  \[
\nabla^{\mu}J_{\mu}(\psi_{1}, \psi_{2})= -\overline{D \psi_{1}}\dual \beta \psi_{2}+ \overline{\psi}_{1}\dual \beta D\psi_{2}, \quad \psi_{i}\in \cinf(M; S(M)).
\]
Using then Stokes formula this implies that 
the Dirac operator $D$ is formally selfadjoint on $\coinf(M; S(M))$ with respect to the indefinite Hermitian form
 \beq\label{e15a.0a}
(\psi_{1}| \psi_{2})_{M}\defeq\int_{M} \overline{\psi}_{1}\dual \beta \psi_{2}\, dV\!\!ol_{g}.
\eeq
\subsubsection{Characteristic manifold}
The {\em principal symbol} $\sigma_{\rm pr}(D)$ equals
\[
\sigma_{\rm pr}(D)(x, \xi)= \gamma(g^{-1}(x)\xi),\ (x, \xi)\in T^{*}M\setminus\zero,
\]
where $\zero= X\times \{0\}$ is the zero section in $T^{*}M$. 

The
 {\em characteristic manifold} of $D$ is
\[
{\rm Char}(D)\defeq\{(x, \xi)\in T^{*}M\setminus\zero: \sigma_{\rm pr}(D)(x, \xi)\hbox{  not invertible}\},
\]
 equal to
\[
{\rm Char}(D)= \{(x, \xi)\in T^{*}M\setminus\zero: \xi\dual g^{-1}(x)\xi=0\}\eqdef \cN,
\]
by the Clifford relations. We denote as usual by $\cN^{\pm}$ the two connected components of $\cN$, where
\beq\label{troup.4}
\cN^{\pm}\defeq\{(x, \xi)\in \cN: \pm \xi\dual v>0\hbox{ for }v\in T_{x}M\hbox{ future directed}\}.
\eeq
\subsubsection{Retarded/advanced inverses}\label{retardo}
Let us assume in addition that $(M, g)$ is globally hyperbolic.
Then, (see \cite{Di} for  Dirac operators in $4$ dimensions, or \cite{M} for more general prenormally hyperbolic operators), $D$ admits unique {\em retarded/advanced inverses} $G_{\rm ret/adv}: \coinf(M; S(M))\to \cinf_{\rm sc}(M; S(M))$ such that: 
\[
\left\{
\begin{array}{l}
DG_{\rm ret/adv}= G_{\rm ret/adv}D= \one, \\[2mm]
 \supp G_{\rm ret/adv}u\subset J_{\pm}(\supp u), \quad u\in \coinf(M; S(M)),
\end{array}\right.
\]
where $J_{\pm}(K)$ are the future/past causal shadows of $K\Subset M$.

 Using the fact that $D$ is formally selfadjoint with respect to $(\cdot| \cdot)_{M}$ and the uniqueness of $G_{\rm ret/adv}$ we obtain that 
\[
G_{\rm ret/adv}^{*}= G_{\rm adv/ret},
\]
where the adjoint is computed with respect to $(\cdot| \cdot)_{M}$. One defines then 
 the {\em causal propagator}
\[
G\defeq G_{\rm ret}- G_{\rm adv}
\] 
which satisfies
\begin{equation}
\label{e15c.0b}\left\{
\begin{array}{l}
DG= GD= 0, \\[2mm]
\supp Gu\subset J(\supp u), \quad u\in \coinf(M; S(M)),\\[2mm]
G^{*}=- G,
\end{array}\right.
\end{equation}
where $J(K)= J_{-}(K)\cup J_{+}(K)$ is the causal shadow of $K\Subset M$.
\subsubsection{The Cauchy problem}\label{sec15c.2}
Let $\Sigma\subset M$ be a smooth, space-like Cauchy surface and denote by $n$ its future directed unit normal and by $S(\Sigma)$ the restriction of the spinor bundle $S(M)$ to $\Sigma$ and
\[
\varrho_{\Sigma}: \cinf(M; S(M))\ni \psi\longmapsto \psi\traa{\Sigma}\in \cinf(\Sigma;S(\Sigma))
\]
the restriction to $\Sigma$. The Cauchy problem
\[
\left\{\begin{array}{l}
D\psi=0,\\
\varrho_{\Sigma}\psi= f, \quad f\in \coinf(\Sigma;S(\Sigma)),
\end{array}\right.
\]
is globally well-posed,  see eg \cite{M}, the solution being denoted by $\psi= U_{\Sigma}f$. We have, see eg \cite[Thm. 19.63]{DG}:
\beq\label{e15c.0c}
U_{\Sigma}f(x)= - \int_{\Sigma}G(x, y)\gamma(n(y))f(y)dV\!\!ol_{h},
\eeq
where $h$ is the Riemannian metric induced by $g$ on $\Sigma$.

We equip $\coinf(\Sigma;S(\Sigma))$ with the indefinite Hermitian form
\beq\label{e15c.2}
(f_{1}| f_{2})_{\Sigma} \defeq \int_{\Sigma}\overline{f}_{1}\dual \beta f_{2}\, dV\!\!ol_{h}. 
\eeq
For $g\in \cE'(\Sigma; S(\Sigma))$, we define $\varrho_{\Sigma}^{*} g\in \cD'(M; S(M))$ by 
\[
\int_{M}\overline{\varrho_{\Sigma}^{*} g}\dual \beta u\, dV\!\!ol_{g}\defeq \int_{\Sigma} \overline{g}\dual \beta \varrho_{\Sigma}udV\!\!ol_{h}, \ u\in\cinf(\Sigma; S(\Sigma)),
\]
i.e. $\varrho_{\Sigma}^{*}$ is the adjoint of $\varrho_{\Sigma}$ with respect to the scalar products $(\cdot|\cdot)_{M}$ and $(\cdot| \cdot)_{\Sigma}$. We can rewrite \eqref{e15c.0c} as
\begin{equation}
\label{e15c.0d}
U_{\Sigma}f= (\varrho_{\Sigma}G)^{*}\gamma(n)f, \quad f\in \coinf(\Sigma; S(\Sigma)).
\end{equation}
 \subsubsection{Cauchy evolution}\label{cauchy-evol}
 Let us assume that $M$ is foliated by a family $(\Sigma_{t})_{t\in \rr}$ of space-like smooth Cauchy surfaces, for example the level sets of a {\em Cauchy time function}, see Subsect. \ref{sec1.1} for the definition. 
 
 Denoting  the restriction of $S(M)$ to $\Sigma_{t}$ by $S(\Sigma_{t})$ and $\varrho_{\Sigma_{t}}$ by $\varrho_{t}$, one can introduce the {\em Cauchy evolution}\[
 U(t,s): \coinf(\Sigma_{s}; S(\Sigma_{s}))\to \coinf(\Sigma_{t}; S(\Sigma_{t})), \ t,s\in \rr
 \]
  defined by
 \[
 U(t,s)f= \varrho_{t} U_{\Sigma_{s}}f
\ f\in \coinf(\Sigma_{s}; S(\Sigma_{s})).
 \]
\subsection{Conformal transformations}\label{sec0.2}
We briefly discuss conformal transformations, and refer to \cite[2.7.2]{St} or \cite[Sect. 17.13]{G} for details. 

Let $c\in \cinf(M)$ with $c(x)>0$ and $\tilde{g}= c^{-2}g$. Then the  spin and spinor bundles for $(M, \tilde{g})$ are identical to those for $(M, g)$.
One has:
\beq\label{troup.5}
\begin{array}{l}
\tilde{\gamma}(X)= c^{-1}\gamma(X), \ \tilde{\beta}= c \beta, \\[2mm]
\tnab^{S}_{C}= \nabla^{S}_{C}- \12 c^{-1}\gamma(X)\gamma (\nabla c)+ \12 c^{-1}X\dual dc ,\\[2mm]
\tilde{\slashed{D}}= c^{\frac{n+1}{2}} \slashed{D} c^{\frac{1-n}{2}},\\[2mm]
\tilde{D}\defeq \tilde{\slashed{D}}+ \tilde{m}=c^{\frac{n+1}{2}} D c^{\frac{1-n}{2}}\hbox{ for }\tilde{m}= cm.
\end{array}
\eeq

\subsection{Quantization of Dirac equation on curved spacetimes}\label{sec0.-2}
We now recall the algebraic quantization of Dirac equations, due to Dimock \cite{Di}.
\subsubsection{$CAR$ $*$-algebras}
 Let $(\cY, \nu)$ be a pre-Hilbert space. The {\em CAR} $*$-{\em algebra} over $(\cY, \nu)$, denoted by ${\rm CAR} (\cY, \nu)$, is the unital complex $*$-algebra generated by elements $\psi(y)$, $\psi^{*}(y)$, $y\in \cY$, with the relations
\begin{equation}
\label{e15.2}
\begin{array}{l}
\psi(y_{1}+ \lambda y_{2})= \psi(y_{1})+ \overline{\lambda}\psi(y_{2}),\\[2mm]
\psi^{*}(y_{1}+ \lambda y_{2})= \psi(y_{1})+ \lambda\psi^{*}(y_{2}), \quad y_{1}, y_{2}\in \cY, \lambda \in \cc, \\[2mm]
[\psi(y_{1}), \psi(y_{2})]_{+}= [\psi^{*}(y_{1}), \psi^{*}(y_{2})]_{+}=0, \\[2mm]
 [\psi(y_{1}), \psi^{*}(y_{2})]_{+}= \overline{y}_{1}\cdot \nu y_{2}\one, \quad y_{1}, y_{2}\in \cY,\\[2mm]
 \psi(y)^{*}= \psi^{*}(y),\quad y\in \cY, 
\end{array}
\end{equation}
where $[A, B]_{+}= AB+ BA$ is the anti-commutator.
\subsubsection{Quasi-free states}
As usual a {\em state} on ${\rm CAR} (\cY, \nu)$ is a linear map $\omega: {\rm CAR} (\cY, \nu)\to \cc$ which is positive and normalized, ie
\[
\omega(A^{*}A)\geq 0, \ \omega(\one)=1, \ A\in {\rm CAR} (\cY, \nu).
\]
-- a state $\omega$ is {\em quasi-free} if:
  \[
 \begin{array}{l}
\omega( \prod_{i=1}^{n}\psi^{*}(y_{i})\prod_{j=1}^{m}\psi(y'_{j}))=0, \hbox{\,\, if }n\neq m,\\[2mm]
 \label{e153.14}\omega( \prod_{i=1}^{n}\psi^{*}(y_{i})\prod_{j=1}^{n}\psi(y'_{j}))= \sum_{\sigma\in S_{n}}{\rm sgn}(\sigma)\prod_{i=1}^{n}\omega(\psi^{*}(y_{i}\psi(y_{\sigma(i)})),
 \end{array}
 \]
 where $S_{n}$ is the set of permutations of $\{1, \dots, n\}$.

-- a quasi-free state is  uniquely determined  by its {\em covariances} $\lambda^{\pm}\in L_{\rm h}(\cY, \cY^{*})$, defined by 
\[
\omega(\psi(y_{1})\psi^{*}(y_{2}))\eqdef \overline{y}_{1}\dual \lambda^{+}y_{2},\quad \omega(\psi^{*}(y_{2})\psi(y_{1}))\eqdef \overline{y}_{1}\dual \lambda^{-}y_{2}, \quad y_{1}, y_{2}\in \cY.
\]
The following two results are well-known, see eg \cite[Sect. 17.2.2]{DG}.
\begin{proposition}
Let $\lambda^{\pm} \in L_{\rm h}(\cY,\cY^*)$. Then the following statements are equivalent :
\begin{enumerate}
\item $\lambda^{\pm}$ are the covariances of a gauge invariant quasi-free state on $\CAR(\cY, \nu)$;
\item $\lambda^{\pm} \geq 0$ and $\lambda^+ + \lambda^- = \nu$.
\end{enumerate}
\end{proposition}

\begin{proposition}\label{hurlub}
A quasi-free state $\omega$ on $\CAR(\cY, \nu)$ is pure if and only if there exist projections $c^{\pm} \in L(\cY)$ such that
\[\lambda^{\pm} = \nu \circ c^{\pm}, c^++c^-= \one.
\]
\end{proposition}
\subsubsection{Pre-Hilbert spaces}\label{pre-hilbert}

We now recall several equivalent pre-Hilbert spaces appearing in the quantization of the Dirac equation.

Let us  denote by  $\Sol(D)$ the space of smooth, space compact solutions of the Dirac equation
\[
D\psi=0.
\]
For $\psi_{1}, \psi_{2}\in \Sol(D)$ we set 
\beq\label{e15.5b}
\overline{\psi}_{1}\dual \nu\psi_{2} \defeq \int_{\Sigma} \i J_{\mu}(\psi_{1}, \psi_{2})n^{\mu} dV\!\!ol_{h}= (\varrho_{\Sigma}\psi_{1}|\i \gamma(n)\varrho_{\Sigma}\psi_{2})_{\Sigma},
\eeq
where $\Sigma$ is a smooth space-like Cauchy surface.

Using that  $\nabla^{\mu}J_{\mu}(\psi_{1}, \psi_{2})=0$ the rhs in \eqref{e15.5b}  is independent on the choice of $\Sigma$. Moreover by \eqref{troup.0} $\nu$ is a positive definite scalar product on $\Sol(D)$.

Setting:
\beq\label{troup.3}
\overline{f}_{1}\dual \nu_{\Sigma}f_{2}\defeq\i \int_{\Sigma}\overline{f}_{1}\dual \beta \gamma(n)f_{2}dV\!\!ol_{h}, 
\eeq
 we obtain that 
\[
\varrho_{\Sigma}: (\Sol(D), \nu)\to (\coinf(\Sigma; S(\Sigma)), \nu_{\Sigma})
\]
is unitary, with inverse $U_{\Sigma}$.

It is also well-known, see eg \cite{Di}, that $G: \coinf(M; S(M))\to \Sol(D)$ is surjective with kernel $D\coinf(M; S(M))$ and  that
\[
G: (\dfrac{\coinf(M; S(M))}{D\coinf(M; S(M))}, \i (\cdot| G\cdot)_{M})\to (\Sol(D), \nu)
\]
is unitary. Summarizing,  the maps
\beq\label{e15c.2b}
\begin{CD}
 (\frac{\coinf(M;S(M))}{D\coinf(M;S(M))}, \i (\cdot| G\cdot)_{M}){\,\overset{G}\longrightarrow\,} (\Sol(D), \nu) {\,\overset{\varrho_\Sigma\,}\longrightarrow\,}(\coinf(\Sigma; S(\Sigma)), \nu_{\Sigma})
\end{CD}
\eeq
are unitary maps between pre-Hilbert spaces.

\subsubsection{$\CAR *$-algebra for Dirac fields}
We denote by ${\rm CAR}(D)$ the $*$-algebra ${\rm CAR}(\cY, \nu)$ for $(\cY, \nu)$ one of the equivalent pre-Hilbert spaces in \eqref{e15c.2b}. 

We use the Hermitian form $(\cdot| \cdot)_{M}$ in \eqref{e15a.0a} to pair $\coinf(M; S(M))$ with $\cD'(M; S(M))$ and to identify continuous sesquilinear forms on $\coinf(M; S(M))$ with continuous linear maps from $\coinf(M; S(M))$ to $\cD'(M; S(M))$.  

We use the Hermitian form $(\cdot| \cdot)_{\Sigma}$ in \eqref{e15c.2} in the same way on the Cauchy surface $\Sigma$.

It is natural to  require a weak continuity of the {\em spacetime covariances} $\Lambda^{\pm}$ of a state $\omega$ on ${\rm CAR}(D)$  defined by:
\[
(u| \Lambda^{+}u)_{M}\defeq \omega(\psi(u)\psi^{*}(u), \ (u| \Lambda^{-}u)\defeq \omega(\psi^{*}(u)\psi(u)), \ u\in \coinf(M; S(M)).
\]
Therefore one considers states on $\CAR(D)$ whose spacetime covariances satisfy:
\begin{equation}
\label{e15.7}
\begin{array}{rl}
{\rm (i)}&\Lambda^{\pm}: \coinf(M; S(M))\to \cD'(M; S(M))\hbox{ are linear continuous},\\[2mm]
{\rm (ii)}& \Lambda^{\pm}\geq 0 \hbox{\,\, with respect to }(\cdot|\cdot)_{M},\\[2mm]
{\rm (iii)}&\Lambda^{+}+ \Lambda^{-}= \i G,\\[2mm]
{\rm (iv)}&D\circ \Lambda^{\pm}= \Lambda^{\pm}\circ D=0.
\end{array}
\end{equation} 
Alternatively, one can define $\omega$ by its {\em Cauchy surface covariances} $\lambda_{\Sigma}^{\pm}$, which satisfy 
\begin{equation}
\label{e15.7b}
\begin{array}{rl}
{\rm (i)}&\lambda^{\pm}_{\Sigma}: \coinf(\Sigma; S(\Sigma))\to \cD'(\Sigma; S(\Sigma))\hbox{ are linear continuous},\\[2mm]
{\rm (ii)}&\lambda^{\pm}_{\Sigma}\geq 0\hbox{ for }(\cdot|\cdot)_{\Sigma},\\[2mm]
{\rm (iii)}&\lambda_{\Sigma}^{+}+ \lambda_{\Sigma}^{-}= \i \gamma(n).
\end{array}
\end{equation}
Using \eqref{e15c.0d} one can show  by the same arguments as for Klein-Gordon fields, see \cite[Prop. 7.5]{GOW} that 
\begin{equation}
\label{e15.7c}
\begin{array}{l}
\Lambda^{\pm}= (\varrho_{\Sigma}G)^{*}\lambda^{\pm}_{\Sigma}(\varrho_{\Sigma}G),\\[2mm]
 \lambda^{\pm}_{\Sigma}= (\varrho_{\Sigma}^{*}\gamma(n))^{*}\Lambda^{\pm}(\varrho_{\Sigma}^{*}\gamma(n)).
\end{array}
\end{equation}
We recall that $S(M)\boxtimes S(M)\xrightarrow{\pi\times \pi}M\times M$  is the vector bundle whose fiber over $(x,x')$ is $End(S_{x}(M), S_{x'}(M))$, with  transition maps naturally  inherited from those of $S(M)$. 

By the Schwartz kernel theorem, we can identify $\Lambda^{\pm}$ with distributional sections in $\cD'(M\times M; S(M)\boxtimes S(M))$, still denoted by $\Lambda^{\pm}$. 
  \subsubsection{The role of the Cauchy evolution}\label{cauchy-evol2}
Recall from \ref{cauchy-evol} that we denoted by $U(t,s)$ the Cauchy evolution associated to a foliation by the Cauchy surfaces $(\Sigma_{t})_{t\in \rr}$.  

If $\omega$ is a quasi-free state on $\CAR(D)$, then denoting by $\lambda^{\pm}(t)$ its Cauchy surface covariances on $\Sigma_{t}$ one has
obviously
\begin{equation}
\label{etroup.-2}
\lambda^{\pm}(t)= U(s,t)^{*}\lambda^{\pm}(s)U(s,t), \ t,s \in \rr.
\end{equation}

 \subsubsection{Hadamard states}
 	
The {\em wavefront set} of $A\in \cD'(M\times M; S(M)\boxtimes S(M))$ is defined in  the natural way: introducing local trivializations of $S(M)$ one can assume that $A\in\cD'(M\times M; M_{N}(\cc))$ where $N= {\rm rank}S(M)$  and  the wavefront set of a matrix valued distribution is simply the union of the wavefront sets of its entries.

We will identify  $T^{*}(M\times M)$ with $T^{*}M\times T^{*}M$. 
If $\Gamma\subset T^{*}M\times T^{*}M$ then one sets
\[
\Gamma'\defeq\{((x, \xi), (x', \xi')): ((x, \xi), (x', -\xi')\in \Gamma\}.
\]
For example $\WF(\delta(x- x'))= \Delta$, where $\Delta\subset T^{*}M\times T^{*}M$ is the diagonal. 

We recall that $\cN^{\pm}$ are the two connected components of $\cN$, see \eqref{troup.4}.

The following definition of Hadamard states is due to Hollands \cite{Ho}.
\begin{definition}\label{def15c.4}
 $\omega$ is a {\em Hadamard state} if 
 \[
\WF(\Lambda^{\pm})'\subset \cN^{\pm}\times \cN^{\pm}.
\]
\end{definition}

The following proposition, see \cite[Prop. 3.8]{St} gives a sufficient condition for the 
Cauchy surface covariances $\lambda^{\pm}_{\Sigma}$ to generate a Hadamard state. 
\begin{proposition}\label{prop15.0a}
 Let \[
 \lambda_{\Sigma}^{\pm}\eqdef \i \gamma(n)c^{\pm}
 \]
  be the Cauchy surface covariances of a quasi-free state $\omega$. Assume that $c^{\pm}$ are  continuous from $\coinf(\Sigma;\cS_{\Sigma})$ to $\cinf(\Sigma;\cS_{\Sigma})$  and from $\cE'(\Sigma; \cS_{\Sigma})$ to $\cD'(\Sigma; \cS_{\Sigma})$, and that for some neighborhood $U$ of $\Sigma$ in $M$ we have
 \beq\label{troup.100}
\WF(U_{\Sigma}\circ c^{\pm})'\subset (\cN^{\pm}\cup\cF)\times T^{*}\Sigma, \hbox{ over }U\times \Sigma,
\eeq
where $\cF\subset T^{*}M$ is a conic set with $\cF\cap \cN= \emptyset$. Then 
$\omega$ is a Hadamard state.
 \end{proposition}

  \subsubsection{Action of conformal transformations}\label{sec5.1b.1}
 Let us now study the action of the conformal transformations recalled in Subsect. \ref{sec0.2}.  
If $\tilde{D}$ is the Dirac operator for $\tilde{g}$, its causal propagator is
\[
\tilde{G}= c^{\frac{n-1}{2}}G c^{-\frac{n+1}{2}}. 
\]
If we set
\[
\begin{array}{l}
W\tilde{\psi}= c^{\frac{1-n}{2}}\tilde{\psi}, \ \tilde{\psi}\in \coinf(M; S(M)),\\[2mm]
W^{*}\psi= c^{\frac{n+1}{2}}\psi, \ \psi\in \coinf(M,; S(M)),\\[2mm]
 Uf= c^{\frac{n-1}{2}}f, f\in  \coinf(\Sigma; S(\Sigma)), 
 \end{array}
\]
then a routine computation gives the following proposition.
\begin{proposition}\label{prop20.1}
 The following diagram is commutative, with all arrows unitary: 
 \[
\begin{CD}
 (\frac{\coinf(M; S(M))}{D\coinf(M; S(M))}, (\cdot\,| \i G\,\cdot)_{M})@>G>> (\Sol(D), \nu) @>
\varrho_{\Sigma}>>(\coinf(\Sigma; S(\Sigma)), \nu_{\Sigma})\\
@VVW^{*}V@VV W^{-1}V@VV  UV\\
 (\frac{\coinf(\tilde{M}; S(M))}{\tilde{D}\coinf(\tilde{M}; S(M))}, (\cdot\,| \i \tilde{G}\,\cdot)_{\tilde{M}})@>\tilde{G}>> (\Sol(\tilde{D}), \tilde{\nu}) @>
\tilde{\varrho}_{\Sigma}>>(\coinf(\Sigma; S(\Sigma)), \tilde{\nu}_{\Sigma})\\
\end{CD}
\]
\end{proposition}
Let us now consider the action of conformal transformations on quasi-free states.

Let $\Lambda^{\pm}$ be the spacetime covariances of a quasi-free state $\omega$ for $D$. Then 
\begin{equation}
\label{e20.3}
\tilde{\Lambda}^{\pm}= c^{\frac{n-1}{2}}\Lambda^{\pm} c^{-\frac{n+1}{2}}
\end{equation}
are the spacetime covariances of a quasi-free state $\tilde{\omega}$ for $\tilde{D}$, 
and
\[
\tilde{\lambda}_{\Sigma}^{\pm}= (U^{*})^{-1}\lambda_{\Sigma}^{\pm}U^{-1}= c^{\frac{n-1}{2}}\lambda_{\Sigma}^{\pm}c^{\frac{1-n}{2}},
\]
if $\lambda_{\Sigma}^{\pm}$, resp. $\tilde{\lambda}_{\Sigma}^{\pm}$ are the Cauchy surface covariances of $\omega$, resp. $\tilde{\omega}$.

Clearly $\omega$ is a Hadamard state iff $\tilde{\omega}$ is.
\subsection{The vacuum state for Dirac fields on static spacetimes}\label{secvac}
The basic example of a state for Dirac fields is the {\em vacuum state} on static spacetimes. Let us recall its definition, following \cite{DH}.

\subsubsection{Vacuum state associated to a Killing field}\label{allo}
Let $(M, g)$ a globally hyperbolic spacetime with a spin structure. The {\em Lie derivative} of a spinor field is defined as  (see \cite{Kos}) :
\beq\label{kosmano}
\begin{array}{l}
\cL_{X}\psi= \nabla_{X}^{S}\psi+ \frac{1}{8}((\nabla_{a}X)_{b}- (\nabla_{b}X)_{a})\gamma^{a}\gamma^{b}\psi,\\[2mm]
 \psi\in \cinf(M; S(M)), \ X\in \cinf(M; TM).
\end{array}
\eeq
If $X$ is a complete  Killing vector field, and  the mass $m$ in \eqref{defdedir} satisfies $X\dual dm=0$, then $[D, \cL_{X}]=0$, see eg \cite[Appendix A]{GHW}.  It follows that the flow $\phi_{s}$ generated by $\cL_{X}$ preserves $\Sol(D)$ and one can easily show, using \eqref{troup.1} and \eqref{kosmano} that it preserves the Hilbertian scalar product $\nu$.  

It hence defines a unique strongly continuous unitary group $(\e^{\i sH})_{s\in \rr}$ on the completion of $(\Sol(D), \nu)$, whose generator $H$ is,  by Nelson's invariant domain theorem, the closure of $\i^{-1}\cL_{X}$ on $\Sol(D)$. 

If $\Sigma$ is a smooth space-like Cauchy surface, we denote by $H_{\Sigma}$ the corresponding generator on the completion of $(\coinf(\Sigma; S(\Sigma)), \nu_{\Sigma})$.

The following definition is taken from \cite{DH}.
\begin{definition}
 Assume that 
 \begin{equation}
\label{kerokero}
\Ker H_{\Sigma}= \{0\}.
\end{equation}
The {\em vacuum state} $\omega^{\rm vac}$  associated to the complete Killing field $X$ is the quasi-free state defined by the Cauchy surface covariances:
\[
\lambda^{\pm{\rm vac}}\defeq \i \gamma(n)\one_{\rr^{\pm}}(H_{\Sigma}).
\]
\end{definition}
Unlike  the bosonic case, $X$ does not need to be time-like in order to be able to define the associated vacuum state.
\subsubsection{Vacuum state on static spacetimes}
We now discuss the  vacuum state on static spacetimes. We will 
 assume that  $M = \rr\times \Sigma$ is equipped with the static metric $g= -c^{2}(\rx)dt^{2}+ h(\rx)d\rx^{2}$, where $c\in\cinf(\Sigma; \rr)$ with $c(\rx)>0$ and $h$ is a Riemannian metric on $\Sigma$.   We set
\[
\tg= c^{-2}g= -dt^{2}+ \th(\rx)d\rx^{2},
\]
which is ultra static.  The restriction of $S(M)$ to $\Sigma_{t}$ is independent on $t$ and denoted by $S(\Sigma)$, see \cite[Subsect. 7.1]{St}.


We consider a {\em static Dirac operator}
\[
D= \slashed{D}+ m,
\]
where $m\in \cinf(\Sigma,; \rr)$  is {\em independent on }$t$.

The corresponding Dirac operator on $(M, \tg)$ is
\[
\tilde{D}= \tilde{\slashed{D}}+ \tilde{m}, \ \tilde{m}= cm. 
\]

 If $(\tilde{e}_{j})_{1\leq j\leq d}$ is a local orthonormal frame for $\tilde{h}$ and $\tilde{e}_{0}= \p_{t}$, we have setting $\tilde{\gamma_{0}}= \tilde{\gamma}(\tilde{e}_{0})$:
 \[
\tilde{D}= -\tilde{\gamma_{0}}(\p_{t}- \i \tilde{H}_{\Sigma})
 \]
 for
 \beq\label{blito}
 \tilde{H}_{\Sigma}= \i \tilde{\gamma_{0}} (\tilde{\gamma}(\tilde{e}_{j}) \tilde{\nabla}^{S}_{\tilde{e_{j}}}+  \tilde{m})\eqdef \tilde{H}_{0\Sigma}+ \i \tilde{\gamma_{0}}\tilde{m}.
 \eeq
 From \eqref{kosmano}, we obtain that $\cL_{\tilde{e_{0}}}= \tilde{\nabla}^{S}_{\tilde{e_{0}}}= \p_{t}$ and hence the generator  of  the Lie derivative w.r.t. the Killing vector field $\p_{t}$ equals $\tilde{H}_{\Sigma}$
 on $\coinf(\Sigma; S(\Sigma))$. We still denote by $\tilde{H}_{\Sigma}$ its closure  for the Hilbertian scalar product $\tilde{\nu}_{\Sigma}$.

 Let us now  consider the original Dirac operator $D$. Using \eqref{troup.5} one  checks that
 \begin{equation}\label{etroup.13}
 D= -c^{-1}\gamma(e_{0})(\p_{t}-\i H_{\Sigma}),
 \end{equation}
 \beq\label{etroup.03}
 H_{\Sigma}\defeq c^{\frac{1-n}{2}}\tilde{H}_{\Sigma}c^{\frac{n-1}{2}},
  \eeq
 where $e_{a}= c^{-1}\tilde{e}_{a}$.  By Prop. \ref{prop20.1} we know that $H_{\Sigma}$  with domain $c^{\frac{1-n}{2}}\Dom \tilde{H}_{\Sigma}$ is selfadjoint for the scalar product $\nu_{\Sigma}$.   It equals the generator of the unitary group associated to $\cL_{\p_{t}}$ considered in \ref{allo}.

 Applying the discussion in \ref{sec5.1b.1} we can define:
 
\begin{definition}\label{def-vaco1}
 Assume that $\Ker H_{\Sigma}= \{0\}$.  Then the  {\em vacuum state} $\omega^{\rm vac}$ for $D$ is the quasi-free state with Cauchy surface covariances
 \[
 \lambda^{\pm{\rm vac}}= \i \gamma(e_{0})\one_{\rr^{\pm}}(H_{\Sigma}).
 \]
 \end{definition}
  \subsubsection{Massive Dirac operators}\label{energy-gap}
  \begin{definition}\label{def-massive-static}
 The static Dirac operator $D$ is called {\em massive} if
 \begin{equation}
 \label{etroup.12}
0\not\in \sigma(H_{\Sigma}).
\end{equation}
\end{definition}
 If is a standard fact  that  if \eqref{etroup.12} holds, then $\omega^{\rm vac}$ is a Hadamard state, see eg. \cite[Thm. 5.1]{SV2}. Another  proof  is given in \cite[Subsect. 7.1]{St}.  If $0\in \sigma(H_{\Sigma})$ but $\Ker H_{\Sigma}= \{0\}$, then one can encounter infrared problems.

  Let us give a simple sufficient condition for \eqref{etroup.12}. Using the Clifford relations and \eqref{troup.1}  we obtain that
 \[
 \tilde{H}_{\Sigma}^{2}=\tilde{H}_{0\Sigma}^{2}+  \tilde{\gamma}(\tilde{h}^{-1}d\tilde{m})+ \tilde{m}^{2}.
 \]
 Since $A=  \tilde{\gamma}(\tilde{h}^{-1}d\tilde{m})$ is selfadjoint for $\tilde{\nu}_{\Sigma}$ with $A^{2}= d\tilde{m}\dual \tilde{h}^{-1}d \tilde{m}$, we obtain that if
 \beq\label{etroup.9}
 \inf_{\Sigma}\tilde{m}^{2}- d\tilde{m}\dual \tilde{h}^{-1}d\tilde{m}>0
 \eeq
  Then $0\not\in \sigma(\tilde{H})$.    In terms of $c, m$ \eqref{etroup.9} becomes:
  \begin{equation}
 \label{etroup.8}
 \inf_{\Sigma} (c^{2}m^{2}- d(cm)\dual h^{-1}d(cm))>0,
 \end{equation}
 Note that \eqref{etroup.8} holds if $c\equiv 1$ and $m(\rx)\equiv m_{0}\neq 0$.

 \section{Dirac operators  on asymptotically static spacetimes}\init\label{sec1}
 \subsection{ Asymptotically static spacetimes}\label{sec1.1}
We fix an orientable  $d-$dimensional manifold $\Sigma$ equipped with a reference Riemannian metric $k$ such that $(\Sigma, k)$ is of bounded geometry, and consider $M= \rr_{t}\times\Sigma_{\rx}$, setting $x= (t, \rx)$, $n= 1+d$ is even.

 \subsubsection{Bounded geometry}
 Roughly speaking a  Riemannian manifold $(\Sigma, k)$ is of bounded geometry if its radius of injectivity is strictly positive and if the metric and all its derivatives, expressed in normal coordinates at a point $\rx$, satisfy estimates which are {\em uniform} with respect to the point $\rx$.

 The two basic examples are compact Riemannian manifolds and $\rr^{d}$ with the flat metric, but many other non compact Riemannian manifolds are of bounded geometry, like for example asymptotically hyperbolic Riemannian manifolds. 
 
 After fixing a background Riemannian metric, one can define in a canonical way various global spaces, like spaces of bounded tensors, Sobolev spaces, bounded differential operators. 
 
 Roughly speaking an object is bounded, if, when expressed in normal coordinates at a base point $\rx$, the object and all its derivatives satisfy estimates which are uniform with respect to $\rx$.
 
 The main interest for us is that on a Riemannian manifold of bounded geometry one can define a global  pseudodifferential calculus, the {\em Shubin calculus}, which shares several important properties with the pseudodifferential calculus on compact manifolds or the uniform pseudodifferential calculus on $\rr^{d}$.
\subsubsection{Lorentzian metric}
We equip $M$ with a Lorentzian metric $g$ of the form
\begin{equation}
\label{e10.1}
g= - c^{2}(x)dt^{2}+ (d\rx^{i}+ b^{i}(x)dt)h_{ij}(x)(d\rx^{j}+ b^{j}(x)dt),
\end{equation} 
where $c\in \cinf(M; \rr),c(x)>0$,  $b\in \cinf(M; T\Sigma)$ and $h\in \cinf(M; \otimes^{2}_{\rm s}T^{*}\Sigma)$ is a  $t$-dependent Riemannian metric on $\Sigma$.

We recall that  $\tilde{t}\in \cinf(M;\rr)$ is called a  \emph{time function} if  $\nabla \tilde{t}$ is a timelike vector field. It is called a \emph{Cauchy time function} if in addition its level sets are Cauchy hypersurfaces.

By \cite[Thm. 2.1]{CC} we know that $(M, g)$ is globally hyperbolic and  $t$ is  a Cauchy time function.

\subsubsection{Asymptotically static spacetimes}\label{sec1.1.1}
We consider also two static metrics  on $M$:
\[
g_{\oi}= -c^{2}_{\oi}(\rx)dt^{2}+ h_{\oi}(\rx)d\rx^{2},
\] 
where $h_{\oi}$, resp. $c_{\oi}$ are two Riemannian metrics, resp. smooth functions  on $\Sigma$ such that:
\[
{\rm (H1)}\ \begin{array}{l}
h_{\oi}\in \BT^{0}_{2}(\Sigma, k), \ h_{\oi}^{-1}\in \BT^{2}_{0}(\Sigma, k), \ c_{\pm\infty}, c_{\pm\infty}^{-1}\in\BT^{0}_{0}(\Sigma, k).
\end{array}
\]
Concerning the asymptotic behavior of $g$ when $t\to \pm \infty$ we assume that
\[
{\rm (H2)}\ \ \beal
h(x)- h_{\oi}(\rx)\in \cS^{-\mu}(\rr^{\pm}; \BT^{0}_{2}(\Sigma, k)), \\[2mm]
b(x)\in \cS^{-1-\mu}(\rr; \BT^{1}_{0}(\Sigma, k)), \\[2mm]
c(x)- c_{\oi}(\rx)\in \cS^{-\mu}(\rr^{\pm}; \BT^{0}_{0}(\Sigma, k)),\\[2mm]
\eeal
\]
for some $\mu>0$, where $\BT^{p}_{q}(\Sigma, k)$ is the  Fr\'echet space  of bounded $q,p$-tensors, see eg \cite{Sh} or \cite[Subsect. 4.1]{St},  and  the space $\cS^{\delta}(\rr; \cF)$  for $\cF$ a Fr\'echet space is defined in Subsect. \ref{sec2.2}. 

In other words the metric $g$ is asymptotic to the static metrics $g_{\oi}$ when $t\to \pm \infty$.

For later use we also fix $m\in\cinf(M; \rr)$, representing a variable mass and $m_{\pm\infty}\in \cinf(\Sigma; \rr)$ such that
\[
{\rm (H3)}\  m(x)- m_{\pm\infty}(\rx)\in \cS^{-\mu}(\rr^{\pm}; \BT^{0}_{0}(\Sigma, k)).
\]
 \subsubsection{Orthogonal decomposition}\label{sec1.1.2}
We recall now the well-known orthogonal decomposition of $g$ associated to the  Cauchy time function $t$.
We set
\[
v\defeq  \dfrac{g^{-1}dt}{dt\cdot g^{-1}dt}=\p_{t}+ b^{i}\p_{\rx^{i}},
\]
which using  ${\rm (H1)}, {\rm (H2)}$ is a complete vector field on $M$. Denoting its flow by $\phi_{t}$ we have:
\[
\phi_{t}(0, \ry)= (t, \rx(t,0, \ry)),\ t\in \rr, \ \ry\in \Sigma,
\]
where $\rx(t, s, \cdot)$ is the flow of the time-dependent vector field $b$ on $\Sigma$.
We also set
\beq\label{e10.0}
\chi: \rr\times \Sigma\ni(t, \ry)\mapsto (t, \rx(t, 0, \ry))\in \rr\times \Sigma.
\eeq
The following lemma is proved in \cite[Appendix A.4]{GW}.
Bounded diffeomorphisms on a manifold of bounded geometry  are defined for example in \cite[Def. 3.3]{GW}.
\begin{lemma}\label{l10.1}Assume ${\rm (H1)}, {\rm (H2)}$. Then
\beq\label{e1.2}
\hat{g}\defeq  \chi^{*}g= -  \hat{c}^{2}(t, \ry)dt^{2}+ \hat{h}(t, \ry)d\ry^{2},
 \eeq
 for $\hat{c}\in \cinf(\rr\times M)$, $\hat{h}\in \cinf(\rr; T^{0}_{2}(\Sigma))$.
 Moreover there exist bounded diffeomorphisms $\rx_{\oi}$ of $(\Sigma, k)$ such that if:
\[
\begin{array}{l}
{\hat h} _{\oi}\defeq \rx_{\oi}^{*}h_{\oi},\\[2mm]
  {\hat c}_{\oi}\defeq  \rx_{\oi}^{*}c_{\oi},
\end{array}
\]
then:
\[
\begin{array}{l}
 \hat h_{\oi}\in \BT^{0}_{2}(\Sigma, k), \  \hat h_{\oi}^{-1}\BT^{2}_{0}(\Sigma, k),\\[2mm]
{\hat c}_{\oi},  {\hat c}_{\oi}^{-1}\in \BT^{0}_{0}(\Sigma, k),
\end{array}
\]
and furthermore,
\[
\begin{array}{rl}
 \hat h- \hat  h_{\oi}\in \cS^{- \mu}(\rr^{\pm}, \BT^{0}_{2}(\Sigma, k)), \\[2mm]
{\hat c}- {\hat c}_{\oi}\in \cS^{- \mu}(\rr^{\pm}, \BT^{0}_{0}(\Sigma, k)),\\[2mm]
\chi^{*}m- m_{\pm\infty}\in \cS^{- \mu}(\rr^{\pm}, \BT^{0}_{0}(\Sigma, k)).
 \end{array}
\]
\end{lemma}

After applying the isometry $\chi: (M, \hat{g})\tosim (M, g)$ in Lemma \ref{l10.1},  removing the hats  to simplify notation  and denoting $\ry$ again by $\rx$, we can assume that 
 \[
  g=  -  c^{2}(t, \rx)dt^{2}+ h(t, \rx)d\rx^{2},
 \]
 with 
 \beq\label{hypo1}
 \begin{array}{l}
   h- h_{\oi}\in \cS^{- \mu}(\rr^{\pm}, \BT^{0}_{2}(\Sigma, k)), \\[2mm]
c- c_{\oi}\in \cS^{- \mu}(\rr^{\pm}, \BT^{0}_{0}(\Sigma, k)),\\[2mm]
m- m_{\pm\infty}\in \cS^{- \mu}(\rr^{\pm}, \BT^{0}_{0}(\Sigma, k)),\\[2mm]
 h_{\oi}\in \BT^{0}_{2}(\Sigma, k), \  h_{\oi}^{-1}\BT^{2}_{0}(\Sigma, k),\\[2mm]
 c_{\oi},   c_{\oi}^{-1}\in \BT^{0}_{0}(\Sigma, k).
 \end{array}
 \eeq
 
 \subsubsection{Conformal transformation}\label{sec1.1.4}
We set  \[
\blibli{g}\defeq c^{-2}g= - dt^{2}+ \blibli{h}(t, \rx)d\rx^{2}
\] 
and obtain that
\beq\label{hypo2}
 \begin{array}{l}
   \blibli{h}- \blibli{h}_{\oi}\in \cS^{- \mu}(\rr^{\pm}, \BT^{0}_{2}(\Sigma, k)),\hbox{ with} \\[2mm]
\blibli{h}_{\oi}= c_{\oi}^{-2}h_{\oi}\in \BT^{0}_{2}(\Sigma, k), \  \blibli{h}_{\oi}^{-1}\BT^{2}_{0}(\Sigma, k).
 \end{array}
 \eeq

\subsection{Spin structures}\label{sec1.2}

 Let us assume that $(M, g)$ admits a spin structure $P\Spin(M, g)$.    We denote by $\Cl(M, g)$, $S(M)$ the Clifford and spinor bundles over $(M, g)$. 
 
By  well-known results on conformal transformations of spin structures, see eg \cite[Lemma 5.27]{LM}, \cite{Hi}\cite[2.7.2]{St}
 $(M, \blibli{g})$ also admits a spin structure and the spinor bundle for $(M, \blibli{g})$ is  equal to $S(M)$.

 Before further discussing the spin structure on $(M, g)$ or $(M, \tilde{g})$ we prove a lemma. We set $\tilde{h}_{t}= \tilde{h}(t, \cdot)$.

\begin{lemma}\label{lemma2.1}
Let us fix a bounded atlas $(V_{i}, \psi_{i})_{i\in \nn}$ for $(\Sigma, \blibli{h}_{0})$.

 Let $\cF_{i}= (e_{i, j})_{1\leq j\leq d}$ oriented orthonormal frames for $\blibli{h}_{0}$ over $V_{i}$ such that $e_{i, j}$ for $i\in \nn$, $1\leq j\leq d$  are a bounded family in $\BT^{1}_{0}(V_{i}, k)$.  Let $\cF_{i}(t)= (e_{i, j}(t))_{1\leq j\leq d}$ the oriented orthonormal frames for $\blibli{h}_{t}$ over $V_{i}$ obtained by parallel transport with respect to $\p_{t}$ of $\cF_{i}$ for the metric $\blibli{g}$. Then:
 
 \ben
 \item $e_{i, j}(\pm\infty)= \lim_{t\to \pm\infty}e_{i, j}(t)$ exist and  the family $e_{i, j}(\pm\infty)$ for $i\in \nn$, $1\leq j\leq d$ is bounded in $\BT^{1}_{0}(V_{i}, k)$.
 
 \item $ \rr^{\pm}\ni t\mapsto e_{i, j}(t)- e_{i, j}(\pm\infty)$ form a bounded family in $\cS^{-\mu}(\rr, \BT^{0}_{0}(V_{i}, k))$.  
 \een
\end{lemma}
\proof 
Let us forget the index $i$ for the moment. Let $x^{\alpha}$, $1\leq \alpha\leq d$ be local coordinates on $V$ obtained from $\psi: V\to B_{d}(0, 1)$ and let $x^{0}= t$. Denoting by $\Gamma^{\mu}_{\varrho\nu}$  the Christoffel symbols  for $g$ in  the local coordinates $(x^{\mu})_{0\leq \mu\leq d}$ over $U= \rr\times V$, we have $\Gamma^{\mu}_{0\nu}= \12 h^{\mu\varrho}\p_{t}h_{\varrho\nu}$. 

Putting back the index $i$ we see from \eqref{hypo2} that $\rr\in t\mapsto \Gamma^{\mu}_{i, 0\nu}(t)$ form a bounded family in $\cS^{-1-\mu}(\rr, \BT^{0}_{0}(V_{i}))$. 
Denoting $e_{i, j}(t)$ simply by $u(t)$ and setting $u= u^{\alpha}\p_{x^{\alpha}}$ over $V$, we obtain that $u(t)$ solves:
\[
\left\{\begin{array}{l}
\p_{t}u^{\alpha}(t)+ \Gamma^{\alpha}_{0\beta}(t)u^{\beta}(t)=0,\\
u^{\alpha}(0)= e_{i, j}^{\alpha}.  
\end{array}
 \right.
\]
From the above estimates on $\Gamma^{\mu}_{i, 0\nu}(t)$ and standard estimates on solutions of linear differential equations,  we obtain (1). It follows that $u(t)$ also solves
\[
\left\{\begin{array}{l}
\p_{t}u^{\alpha}(t)+ \Gamma^{\alpha}_{0\beta}(t)u^{\beta}(t)=0,\\
\lim_{t\to \pm \infty} u^{\alpha}(t)= u^{\alpha}(\pm\infty).  
\end{array}
 \right.
\]
Again the same estimates (integrating now from $t=\pm \infty$ instead of from $t=0$) prove (2) and  complete the proof of the lemma. 
 \qed

 \subsubsection{Spin structures}\label{sec1.2.1}

 Since $M$ is a cartesian product and from the form of $\tilde{g}$, further simplications occur, see eg \cite{BGM} or \cite[Subsect. 2.6]{St}.

 Let us set 
$  \bar{\rr}= \rr\cup\{-\infty, +\infty\}$
 and set $\tilde{h}_{{\rm in/out}}= \tilde{h}_{\mp\infty}$ for coherence of notation.
 
We can use the local frames $\cF_{i}(t)$ over $V_{i}$  to obtain local trivialisations of $P\SO(\Sigma, \tilde{h}_{t})$ for $t\in \bar{\rr}$. The associated transition functions are independent on $t$. By the arguments in \cite[Subsect. 2.6]{St}, we obtain unique spin structures on $(\Sigma, \tilde{h}_{t})$ for $t\in \bar{\rr}$.
 
 The transition functions of $P\Spin(M, \tilde{g})$ are independent on $t$ and induce 
 a spin structure on $(\Sigma, \tilde{h}_{t})$ whose transition functions are also independent on $t$. 
 If $S_{t}(\Sigma)$ denotes the restriction of $S(M)$ to $\Sigma_{t}$, then $S_{t}(\Sigma)$ is independent on $t$ and denoted by $S(\Sigma)$.
 
 Conversely the spin structure on $(\Sigma, \tilde{h}_{\pm\infty})$ induces a spin structure on $(\Sigma, \tilde{g}_{\pm\infty})$ for $\tilde{g}_{\pm\infty}= -dt^{2}+ \tilde{h}_{\pm\infty}(\rx)d\rx^{2}$ and by conformal invariance a spin structure on $(M, g_{\outin})$. The associated spinor bundle is again equal to $S(M)$.

\subsection{Dirac operators}
 We consider the Dirac operator locally given by
  \begin{equation}
  \label{e1.2b}
D\defeq  \slashed{D}+m, \  \slashed{D}=g^{ab}\gamma(e_{a})\nabla^{S}_{e_{b}}
  \end{equation}
  where $(e_{a})_{0\leq  a\leq d}$ is some local frame of $TM$ and  $\nabla^{S}$ is the spin connection. 
  
  \subsubsection{Conformal transformation}\label{confoconfo}
By Subsect. \ref{sec0.2} we obtain that 
\begin{equation}
\label{econfo.1}
D= c^{-\frac{n+1}{2}}\blibli{D}c^{\frac{n-1}{2}} \hbox{ for }\blibli{D}=\slashed{\blibli D}+ \blibli{m}, \ \blibli{m}= cm,
\end{equation}
with
\begin{equation}
\label{hypo3}
\begin{array}{l}
\blibli{m}- \blibli{m}_{\oi}\in  \cS^{- \mu}(\rr^{\pm}, \BT^{0}_{0}(\Sigma, k)), \\[2mm]
 \blibli{m}_{\oi}= c_{\oi}m_{\oi}, \  \blibli{m}_{\oi},  \blibli{m}_{\oi}^{-1}\in \BT^{0}_{0}(\Sigma, k).
\end{array}
\end{equation}
\subsubsection{Asymptotic Dirac operators}
Let
\[
D_{\oi}= \slashed{D}_{\oi}+ m_{\oi}
\]
the asymptotic Dirac operators obtained from the spin structures $P\Spin(M, g_{\oi})$.

 We will  assume
\[
{\rm (H4)}\ D_{\oi}\hbox{ are {\em massive} ie }0\not\in \sigma(H_{\oi}),
\]
see \ref{energy-gap}. A sufficient condition for ${\rm (H4)}$ is given in  \eqref{etroup.8}.

\section{Pseudodifferential calculus}\init\label{sec2}
In this section we will recall Shubin's global pseudodifferential calculus on manifolds of bounded geometry  and its time-dependent versions. We refer the reader to \cite{Sh, Ko} for the original exposition and  to  \cite{GOW} for a more recent one. We are interested in pseudodifferential operators acting on sections of spinor bundles, which are considered in \cite{St}.

\subsection{Notations}\label{sec2.1}
Let $(\Sigma, k)$ a Riemannian manifold of bounded geometry see \cite{CG, Ro} or \cite[Thm. 2.2]{GOW} for an equivalent definition. 
We refer the reader to \cite[Subsect. 4.1]{St} for the definitions below.

We denote by $\BT^{p}_{q}(\Sigma, k)$  the space of bounded $(q, p)$ tensors on $\Sigma$. Let also  $E\xrightarrow{\pi}\Sigma$ a vector bundle of bounded geometry.   

We denote by $S^{m}_{\rm ph}(T^{*}\Sigma; L(E))$ the space of $L(E)$-valued polyhomogenous symbols of order $m$ on $\Sigma$, see eg \cite[Sect. 4.1]{St}.  

The ideal of smoothing operators is denoted by $\cW^{-\infty}(\Sigma; L(E))$, and  one sets
\[
\Psi^{m}(\Sigma; L(E))= \Op(S^{m}_{\rm ph}(T^{*}\Sigma; L(E)))+ \cW^{-\infty}(\Sigma; L(E)),
\]
for some quantization map $\Op$ obtained from a bounded atlas and bounded partition of unity of $(\Sigma, k)$.

\subsection{Time-dependent pseudodifferential operators}\label{sec2.2}
We will also consider time-dependent pseudodifferential operators, adapted to the geometric situation considered in Subsect. \ref{sec1.2}. 

We first introduce some notation.

Let $\mathcal{F}$  a Fr\'echet space whose topology is defined by the seminorms $\| \cdot\|_{p}$, $p\in \nn$ and $\delta\in \rr$. We denote by $\cS^{\delta}(\rr; \mathcal{F})$ the space of smooth functions $f: \rr\to \cF$ such that $\sup_{\rr}\langle t\rangle^{k- \delta}\|\p_{t}^{k}f(t)\|_{p}<\infty$ for all $k, p\in \nn$. Equipped with the obvious seminorms, it is itself a Fr\'echet space.

Note that  $\cS^{\delta}(\rr; \cF)= \langle t\rangle^{\delta}\cS^{0}(\rr; \cF)$ so we can always reduce ourselves to $\delta=0$.

Similarly we denote by  $\cinfb(\rr; \cF)$ the space of smooth functions $f: \rr\to \cF$ such that $\sup_{\rr}\|\p_{t}^{k}f(t)\|_{p}<\infty$ for all $k, p\in \nn$, with the analogous Fr\'echet space topology.

We use this notation to define the spaces $\cS^{\delta}(\rr; S^{m}_{\rm ph}(T^{*}\Sigma; L(E)))$, \hfill\linebreak $\cS^{\delta}(\rr, \cW^{-\infty}(\Sigma; L(E)))$ and $\cS^{\delta}(\rr; \Psi^{m}(\Sigma; L(E)))$.

 For example if $(\Sigma, k)$ equals $\rr^{n}$ equipped with the flat metric, then \hfill\linebreak $\cS^{\delta}(\rr; S^{m}_{\rm ph}(T^{*}\rr^{n}))$ is the space of smooth functions $a: \rr\times T^{*}\rr^{n}\to \cc$ such that  there exist  for $j\in \nn$ functions $a_{m-j}: \rr\times T^{*}(\rr^{n})\to \cc$, homogeneous of degree $m-j$ in $\xi$ with 
\[
\sup_{\rr\times T^{*}\rr^{n}\setminus\zero}\langle t\rangle^{-\delta+ k}\langle \xi\rangle^{-m+j|\beta|}|\p_{t}^{k}\p_{x}^{\alpha}\p_{\xi}^{\beta}a_{m-j}(t, x, \xi)|<\infty, \ k\in\nn, \alpha, \beta\in \nn^{n},
\]
and for any $N\in \nn$
\[
\sup_{\rr\times T^{*}\rr^{n}\setminus\zero}\langle t\rangle^{-\delta+ k}\langle \xi\rangle^{-m+N+1+ |\beta|}|\p_{t}^{k}\p_{x}^{\alpha}\p_{\xi}^{\beta}(a-\sum_{j=0}^{N}a_{m-j}(t, x, \xi))|<\infty, \ k\in\nn, \alpha, \beta\in \nn^{n}.
\]
Similarly  $\cS^{\delta}(\rr, \cW^{-\infty}(\rr^{n}))$ is the space of smooth functions $a: \rr\to B(L^{2}(\rr^{n}))$ such that
\[
\sup_{\rr}\langle t\rangle^{-\delta+ k}\| \p_{t}^{k}a(t)\|_{B(H^{-m}(\rr^{n}), H^{m}(\rr^{n}))}<\infty, \ k, m\in \nn,
\]
where $H^{m}(\rr^{n})$ are the usual Sobolev spaces.

For simplicity of notation  $\cS^{\delta}(\rr; S^{m}_{\rm ph}(T^{*}\Sigma; L(E)))$ or $\cS^{\delta}(\rr; \Psi^{m}(\Sigma; L(E)))$ will often simply be denoted by $\cS^{\delta, m}$, $\Psi^{\delta, m}$.

\subsubsection{Principal symbol}\label{sec2.2.1}
If $A(t)= \Op(a(t))+ R_{-\infty}(t)\in \cS^{\delta}(\rr; \Psi^{m}(\Sigma; L(E)))$  its principal symbol  is
 \[
 \sigma_{\rm pr}(A)(t)\defeq [a](t)\in\cS^{\delta}(\rr;  S^{m}_{\rm ph}(T^{*}\Sigma; L(E)))/\cS^{\delta}(\rr;  S^{m-1}_{\rm ph}(T^{*}\Sigma; L(E))).
 \] 
  $\sigma_{\rm pr}(A)(t)$ is independent on the decomposition of $A(t)$ as $\Op(a)(t)+R_{-\infty}(t)$ and on the choice of the good quantization map $\Op$.
  As usual we choose a representative of  $ \sigma_{\rm pr}(A)(t)$ which is homogeneous of order $m$ on the fibers of $T^{*}\Sigma$.
\subsubsection{Ellipticity}\label{sec2.2.2}
 An operator $A(t)\in \cS^{\delta}(\rr; \Psi^{m}(\Sigma; L(E)))$ is {\em elliptic} if $\sigma_{\rm pr}(A)(t, x, \xi)$ is invertible for all $t\in \rr$ and 
 \[
 \sup_{t\in \rr, (x, \xi)\in T^{*}\Sigma,  | \xi|= 1}\| \sigma_{\rm pr}(A)^{-1}(t, x, \xi)\| <\infty.
 \]
 To define the norm above, one chooses a bounded Hilbert space structure on the fibers of $E$, the definition being independent on its choice.
 
  \begin{proposition}\label{prop4.1time}
Let $A(t)\in \cS^{\epsilon}(\rr;\Psi^{m}(\Sigma; L(E)))$, $\epsilon\in\rr, m\geq 0$ elliptic. Then the following holds:
\ben
\item  $A(t)$ is closeable on $\coinf(\Sigma; E)$ with $\Dom A^{\rm cl}(t)= H^{m}(\Sigma; L(E))$.  
\item  if  there exists $\delta>0$ such that $[-\delta, \delta]\cap  \sigma(A^{\rm cl}(t))=\emptyset$ for $t\in \rr$, then $A^{-1}(t)\in \cS^{-\epsilon}(\rr;\Psi^{-m}(\Sigma; L(E)))$ and
\[
\sigma_{\rm pr}(A^{-1})(t)= (\sigma_{\rm pr}(A))^{-1}(t).
\]
\een
\end{proposition}
\proof the same result is proved in \cite[Prop. 5.8]{St}, with $\cS^{\delta}(\rr; \Psi^{m})$ replaced by $\cinfb(\rr; \Psi^{m})$, where $\cinfb(\rr; \cF)$ is  defined at the beginning of Subsect. \ref{sec2.2}. Note that $a(t)\in \cS^{\delta}(\rr; \cF)$ iff $\langle t\rangle^{-\delta- n}\p_{t}^{n}a(t)\in \cinfb(\rr; \cF)$ for all $n\in \nn$.  Using that  $\p_{t}A^{-1}(t)= - A^{-1}(t)\p_{t}A(t)A^{-1}(t)$ and similar identities for higher derivatives of $A^{-1}(t)$ combined with the above remark we obtain the proposition. \qed
 \subsection{Functional calculus}\label{sec2.3}
  \subsubsection{Elliptic selfadjoint operators}\label{sec4.2.5}
 Let us fix a bounded Hilbertian structure $(\cdot| \cdot)_{E}$ on the fibers of $E$ and  define the scalar product
 \[
 (u|v)= \int_{\Sigma}(u(x)| v(x))_{E}dVol_{g}, \ u, v\in \coinf(\Sigma; E).
 \]
Let  $H(t)\in \cS^{\delta}(\rr;\Psi^{m}(\Sigma; L(E)))$ be  elliptic, symmetric on $\coinf(\Sigma; E)$. Using Prop. \ref{prop4.1time} one easily shows that  its closure is selfadjoint with domain $H^{m}(\Sigma; E)$. Note also that its principal symbol $\sigma_{\rm pr}(H)(t, x, \xi)$ is selfadjoint for the Hilbertian scalar product on $E_{x}$. 
 \subsubsection{Functional calculus}
 
 We now extend some results in \cite{St} on functional calculus for selfadjoint pseudodifferential operators to our situation. We first recall some definitions from \cite[Subsect. 5.3]{St}  about  pseudodifferential operators with parameters.
 
  One denotes by $\tilde{S}^{m}(\Sigma; L(E))$  the space of symbols $b\in\cinf(\rr_{\lambda}\times T^{*}\Sigma; L(E))$  such that if $b(\lambda)= b(\lambda, \cdot)\in \cinf(T^{*}\Sigma; L(E))$ and $T_{i}b(\lambda)$ are the  pushforwards of $b(\lambda)$ associated to a covering $\{U_{i}\}_{i\in \nn}$ of $\Sigma$, we have: 
  \[
\p^{\gamma}_{\lambda}\p^{\alpha}_{x}\p^{\beta}_{\xi}b_{i}(\lambda, x, \xi)\in O(\langle \xi\rangle + \langle \lambda\rangle)^{m -|\beta|- \gamma},\ (\lambda, x, \xi)\in \rr\times T^{*}B(0, 1)
  \]
  uniformly with respect to $i\in \nn$.   One denotes by $\widetilde{S}^{m}_{\rm h}(T^{*}\Sigma; L(E))$ the subspace of such symbols which are homogeneous w.r.t. $(\lambda, \xi)$ and by $\widetilde{S}^{m}_{\rm ph}(T^{*}\Sigma; L(E))$ the subspace of polyhomogeneous symbols.
  
   One also defines the ideal  $\widetilde{\cW}^{-\infty}(\Sigma; L(E))$ as the set of smooth functions $b: \rr\in \lambda\mapsto b(\lambda)\in \cW^{-\infty}(\Sigma; L(E))$ such that
  \[
  \| \p^{\gamma}_{\lambda}b(\lambda)\|_{B(H^{-m}(\Sigma), H^{m}(\Sigma))}\in O(\langle \lambda\rangle^{-n}), \ \forall, m,n,\gamma\in \nn,
  \]
  and set
  \[
  \widetilde{\Psi}^{m}(\Sigma; L(E))\defeq \Op (\widetilde{S}_{\rm ph}^{m}(T^{*}\Sigma; L(E)))+ \widetilde{\cW}^{-\infty}(\Sigma; L(E)).
  \]
  
  As usual one defines  the time-dependent versions of the above spaces:
  \[
  \cS^{\delta}(\rr; \widetilde{S}^{m}_{\rm ph}(T^{*}\Sigma; L(E))), \ \cS^{\delta}(\rr; \widetilde{\cW}^{-\infty}(\Sigma; L(E))), \cS^{\delta}(\rr; \widetilde{\Psi}^{m}(\Sigma; L(E))).
  \]
  We define the {\em principal symbol} of  $A(t)\in \cS^{\delta}(\rr; \widetilde{\Psi}^{m}(\Sigma; L(E)))$ as in \ref{sec2.2.1},  using the polyhomogeneity.
  \begin{proposition}\label{propopi}
 Let $H(t)\in \cS^{\delta}(\rr; \Psi^{1}(\Sigma; L(E)))$ elliptic and formally selfadjoint.  Let us still denote by $H(t)$ its closure, which is selfadjoint on $H^{1}(\Sigma; E)$ by  Prop. \ref{prop4.1time}. Assume that there exists $\delta>0$ such that $[-\delta, \delta]\cap \sigma(H(t))= \emptyset$ for $t\in I$.
 
  Then $\lambda\mapsto (H(t)+ \i \lambda)^{-1}$ belongs to $\cS^{\delta}(\rr; \widetilde{\Psi}^{-1}(\Sigma; L(E)))$ with principal symbol $(\sigma_{\rm pr}(H(t))+ \i \lambda)^{-1}$.
\end{proposition}
\proof The $\cinfb$ version of the proposition  is proved in \cite[Prop. 5.9]{St}.  We use the same remark as in the proof of Prop. \ref{prop4.1time} to extend it to the $\cS^{\delta}$ case. Details are left to the reader. \qed

  \begin{proposition}\label{prop4.3}
Let  $H(t)\in \cS^{0}(\rr; \Psi^{1}(\Sigma; L(E)))$ be elliptic, symmetric on $\coinf(\Sigma; E)$, and let us denote still by $H(t)$ its selfadjoint closure.  Assume that there exists $\delta>0$ such that $[-\delta, \delta]\cap \sigma(H(t))= \emptyset$ for $t\in \rr$. 

Assume in addition that there exist $H_{\infty}\in \Psi^{1}(\Sigma; L(E))$, elliptic symmetric on $\coinf(\Sigma; E)$ with $0\not\in \sigma(H_{\infty})$ such that
\[
H(t)- H_{\infty}\in \cS^{-\mu}(\rr; \Psi^{1}(\Sigma; L(E))).
\]
 Then
\ben
\item  the spectral projections $\one_{\rr^{\pm}}(H(t))$ belong to $\cS^{0}(\rr;\Psi^{0}(\Sigma; L(E)))$ and
\[
\sigma_{\rm pr}(\one_{\rr^{\pm}}(H(t)))= \one_{\rr^{\pm}}(\sigma_{\rm pr}(H(t))).
\]
Moreover $\one_{\rr^{\pm}}(H(t))-\one_{\rr^{\pm}}(H_{\infty})$ belongs to $\cS^{-\mu}(\rr;\Psi^{0}(\Sigma; L(E)))$. 
\item $S(t)= (H^{2}(t)+1)^{\12}$ belongs to  $\cS^{0}(\rr; \Psi^{1}(\Sigma; L(E)))$  and 
\[
\sigma_{\rm pr}(S(t))=|\sigma_{\rm pr}(H(t))|.
\]
Moreover $S(t)- S_{\infty}$ belongs to $\cS^{-\mu}(\rr;\Psi^{1}(\Sigma; L(E)))$ for $S_{\infty}= (H^{2}_{\infty}+1)^{\12}$.
\een
\end{proposition}
\proof 
By Prop. \ref{propopi} we have 
 \beq\label{e4.3}
 (\i \lambda- H(t))^{-1}= \Op(a(t, \lambda))+ R_{-\infty}(t, \lambda),
 \eeq
 where $a(t)\in  \cS^{0}(\rr; \widetilde{S}^{-1}(T^{*}\Sigma; L(E)))$   and $R_{-\infty}(t)\in \cS^{0}(\rr; \widetilde{\cW}^{-\infty}(\Sigma; L(E))$ satisfies:
 \[
\langle t\rangle^{ p}\| \p_{\lambda}^{n}\p_{t}^{p}R_{-\infty}(t, \lambda)\|_{B(H^{-m}(\Sigma), H^{m}(\Sigma))}\in O(\langle \lambda\rangle)^{-m}, \ \forall p, m,n\in \nn,
 \]
 uniformly for $t\in \rr$.
 
  The principal symbol of $a(t)$ is $(\i \lambda- \sigma_{\rm pr}(H))^{-1}$, which means that   \beq\label{e4.4}
  \Op(a(t, \lambda))- \Op((\i \lambda- \sigma_{\rm pr}(H)(t))^{-1}\in \cS^{0}(\rr; \widetilde{\Psi}^{-2}(\Sigma; L(E))).
  \eeq
  For $a\neq 0$ we have
\beq\label{tuod}
|a|^{-1}= \frac{2}{\pi}\int_{0}^{+\infty}(a+ \i\lambda)^{-1}(a-\i \lambda)^{-1}d\lambda,
\eeq
hence
\begin{equation}
\label{tud}
|H(t)|^{-1}=  \frac{2}{\pi}\int_{0}^{+\infty}(H(t)+ \i \lambda)^{-1}(H(t)- \i \lambda)^{-1}d\lambda.
\end{equation}
 From Prop. \ref{propopi}  we obtain that $|H(t)|^{-1}\in \cS^{0}(\rr; \Psi^{-1}(\Sigma; L(E)))$. We also deduce from \eqref{tud} using the second resolvent formula that $|H(t)|^{-1}- |H_{\infty}|^{-1}\in \cS^{-\mu}(\rr; \Psi^{-1}(\Sigma; L(E)))$.  This implies that ${\rm sgn}(H(t))\in \cS^{0}(\rr; \Psi^{-0}(\Sigma; L(E)))$ and ${\rm sgn}(H(t))-{\rm sgn}(H_{\infty})\in \cS^{-\mu}(\rr; \Psi^{0}(\Sigma; L(E)))$.  
 
 Moreover since the principal symbol of $(H(t)+\i\lambda)^{-1}$ equals $(\sigma_{\rm pr}(H(t))+ \i \lambda)^{-1}$, applying once more \eqref{tuod} we obtain that 
 $\sigma_{\rm pr}({\rm sgn}(H(t)))$ equals ${\rm sgn}(\sigma_{\rm pr}(H(t)))$.

 Writing $\one_{\rr^{\pm}}(\lambda)= \12 (1\pm {\rm sgn}(\lambda))$ this implies (1). To prove (2) we deduce from \eqref{tuod} that
 \beq\label{tuad}
   (a+1)^{-\12}= \frac{2}{\pi}\int_{0}^{+\infty}(a+s^{2}+1)^{-1}ds= \frac{2}{\pi}\int_{1}^{+\infty}(a+ \lambda^{2})^{-1}\lambda(\lambda^{2}-1)^{-\12}d\lambda,
\eeq
 hence
\beq\label{tut}
  (H^{2}(t)+1)^{-\12}= \frac{2}{\pi}\int_{1}^{+\infty}(H(t)+ \i \lambda)^{-1}(H(t)- \i \lambda)^{-1}\lambda (\lambda^{2}-1)^{-\12}d\lambda.
  \eeq
 we obtain that $(H^{2}(t)+ 1)^{-\12}\in \cS^{0}(\rr; \Psi^{-1}(\Sigma; L(E)))$. We also deduce from \eqref{tut} that $(H^{2}(t)+ 1)^{-\12}- (H^{2}_{\infty}+ 1)^{-\12}\in \cS^{-\mu}(\rr; \Psi^{-1}(\Sigma; L(E)))$.
  We write then $(H^{2}(t)+1)^{\12}= (H^{2}+1)(H^{2}(t)+1)^{-\12}$ and obtain (2). \qed

    \section{The in/out vacuum states}\init\label{sec4}
 
In this section we prove Thm. \ref{mainmain}.
 \subsection{Reduction of the Dirac operator}\label{sec4.1}
 In this subsection we consider the Dirac operator $\blibli{D}$ obtained from $D$ by conformal transformation, see \ref{sec1.1.4}, \ref{confoconfo}.
 In order not to overburden the notation, we remove the tildes and denote $\blibli{g}$, $\blibli{D}$ etc simply by $g$, $D$.

 We recall that   the restriction $S_{t}(\Sigma)$  of the spinor bundle $S(M)$ to $\Sigma_{t}$ is independent of $t$, and denoted by $S(\Sigma)$.

We recall also that $S(M)$ is equipped by a time positive Hermitian form $\beta$  see 
\eqref{troup.0} and we denote by $\beta_{t}$ its  restriction  to $S(\Sigma)$. Also  we denote by $\gamma_{t}: T_{\Sigma_{t}}M\to L(S_{t}(\Sigma))$ the restrictions of $\gamma$ to $S(\Sigma_{t})$.

We will denote by $(\rx, k)$ local coordinates on $T^{*}\Sigma$ and by $(t, \rx, \tau,k)$ local coordinates on $T^{*}M$.

The first step consists in reducing the Dirac equation $\blibli{D}\psi=0$ to a  time-dependent Schroedinger equation
\[
\p_{t}\psi- \i \mathsf{H}(t)\psi=0,
\]
where $\mathsf{H}(t)$ is some time-dependent selfadjoint operator. To this end it is necessary to identify the elements of spinor bundles at different times by parallel transport.  We  recall that  $e_{0}= \p_{t}$ and $(e_{j})_{1\leq j\leq d}$ are the local frames  constructed in Lemma \ref{lemma2.1}. 
We start by an easy proposition.

\begin{proposition}\label{prop4.0}
\[
\begin{array}{rl}
i)&\gamma_{t}(e_{0})- \gamma_{\pm\infty}(e_{0})\in \cS^{-\mu}(\rr^{\pm}; \cinfb(V; L(S(\Sigma)))),\\[2mm]
ii)&\gamma_{t}(e_{j}(t))- \gamma_{\pm\infty}(e_{j}(\pm\infty))\in \cS^{-\mu}(\rr^{\pm}; \cinfb(V; L(S(\Sigma)))),\\[2mm]
iii)&\beta_{t}- \beta_{\pm\infty}\in \cS^{-\mu}(\rr^{\pm}; \cinfb(V; L(S(\Sigma), S(\Sigma)^{*}))),\\[2mm]
iv)&\nabla^{S}_{e_{0}}- \nabla^{S_{\pm\infty}}_{e_{0}}\in \cS^{-\mu}(\rr^{\pm}; \Psi^{1}(\Sigma; S(\Sigma))),\\[2mm]
v)&\nabla^{S}_{e_{j}(t)}- \nabla^{S_{\pm\infty}}_{e_{j}(\pm\infty)}\in \cS^{-\mu}(\rr^{\pm}; \Psi^{1}(\Sigma; S(\Sigma))).
\end{array}
\]
\end{proposition}
\proof
We fix a bounded atlas $(V_{i}, \psi_{i})_{i\in \nn}$ of $(\Sigma, h_{0})$ and set $U_{i}= \rr\times V_{i}$.   We fix a bounded family $(\cF_{i})_{i\in \nn}$ of oriented orthonormal frames for $h_{0}$ over $V_{i}$ and denote by $\cF_{i}(t)= (e_{i, j}(t))_{1\leq j\leq d}$ the orthonormal frames obtained by parallel transport as in Lemma \ref{lemma2.1}.   Since  $e_{0}= \p_{t}$,  $\cE_{i}= (e_{i, a})_{0\leq a\leq d}$ are then oriented, time-oriented orthonormal frames for $g$ over $U_{i}= \rr\times V_{i}$.

We use  the spin frames $\cB_{i}(t)= (E_{i, A}(t))_{1\leq A\leq N}$ of $S(\Sigma)$ associated to the frames $\cE_{i}(t)= (e_{i, a}(t))_{0\leq a\leq d}$ over $\{t\}\times V_{i}$.  From the estimates  in Lemma \ref{lemma2.1}, we obtain that $\cB_{i}(\pm\infty)= \lim_{t\to \pm \infty}\cB_{i}(t)$ exist and that
\beq\label{et.2}
\begin{array}{l}
E_{i, A}(\pm\infty)\in \cinfb(V_{i}, S(\Sigma)), \\[2mm]
E_{i, A}(t)- E_{i, A}(\pm\infty)\in \cS^{-\mu}(\rr^{\pm}; \cinfb(V_{i}, S(\Sigma)))\hbox{ uniformly w.r.t. }i\in \nn.
\end{array}
\eeq
 We recall from \ref{sec1.2.1} that  from the transition functions $\pmb{o}_{ij}(\pm\infty): V_{ij}\to \SO(d)$    one obtains the spin structures $P\Spin(\Sigma; h_{\pm\infty})$ introduced above. The frames $\cB_{i}(\pm\infty)$ are the frames associated to the $\cE_{i}(\pm\infty)$ for this spin structure.

 Let us  now forget the index $i$ and  denote by $(\pmb{\psi}^{A})_{1\leq A\leq N}\in \cc^{N}$ the components of $\psi$ in the frame $\cB$. 
 The dual frames are as usual denoted by $(e^{a})_{0\leq a\leq d}$, $(E^{A})_{1\leq A\leq N}$ so for example $\pmb{\psi}^{A}= \psi\dual E^{A}$.
 
 Denoting by $\pmb{\gamma}_{t}(u)$ the matrix of $\gamma_{t}(u)$ in the frame $\cB(t)$, we have also
\beq\label{et.2c}
\pmb{\gamma}_{t}(u)= \pmb{\gamma}_{a}u^{a}(t),\ u^{a}(t)\defeq u\dual e^{a}(t).
\eeq
where $\pmb{\gamma}_{a}\in M_{N}(\cc)$ for $0\leq a\leq d$ are the usual gamma matrices.  Using Lemma \ref{lemma2.1} and \eqref{et.2} we obtain that $\lim_{t\to\pm\infty}\gamma_{t}(e_{a}(t))\in L(S(\Sigma))$ exist and that
\begin{equation}
\label{et.2b}
\gamma_{t}(e_{a}(t))- \lim_{t\to\pm\infty}\gamma_{t}(e_{a}(t))\in \cS^{-\mu}(\rr^{\pm}; \cinfb(V; L(S(\Sigma)))).
\end{equation}
 If we reintroduce the index $i$ and set $V= V_{i}$, then the seminorms in \eqref{et.2b} are uniform with respect to $i\in \nn$.  Because of \eqref{et.2c}, the limits  $\lim_{t\to\pm\infty}\gamma_{t}(e_{a}(t))$ are equal to $\gamma_{\pm\infty}(e_{a}(\pm\infty))$. This proves {\it ii)}. {\it i)} is proved similarly.

Let us now denote  by $\pmb{\beta}_{t}$ the matrix of $\beta_{t}$ in the frame $\cB(t)$. We have
\[
\pmb{\beta}_{t}= \pmb{\beta}
\]
where $\pmb{\beta}\in M_{N}(\cc)$ is  a Hermitian matrix such that  
\[
\pmb{\beta}\pmb{\gamma}_{a}= - \pmb{\gamma}_{a}^{*}\pmb{\beta}, \ \i \pmb{\beta}\pmb{\gamma}_{0}>0.
\]
This implies as above that $\lim_{t\to \pm \infty}\beta_{t}$ exist and
\begin{equation}
\label{et2c}
\beta_{t}- \lim_{t\to \pm \infty}\beta_{t}\in \cS^{-\mu}(\rr^{\pm}; \cinfb(V; L(S(\Sigma), S(\Sigma)^{*}))).
\end{equation}
Again because of \eqref{et.2c} the limits $\lim_{t\to \pm \infty}\beta_{t}$ are equal to $\beta_{\pm\infty}$, which proves {\it iii)}.

Finally, see eg \cite[2.5.6]{St}   we have:
\beq\label{et.3}
\nabla^{S}_{e_{a}}\pmb{\psi}^{A}=  \p_{a}\pmb{\psi}^{A}+ \sigma^{A}_{aC}\pmb{\psi}^{C},
\eeq
where
\[
\p_{a}f= e_{a}\dual df, \ \sigma^{A}_{aC}= E^{A}\dual \sigma_{a}E_{C}, \ \sigma_{a}= \frac{1}{4} \Gamma^{c}_{ab}\gamma(e_{c})g^{bd}\gamma(e_{d}),\ \Gamma^{c}_{ab}= \nabla_{e_{a}}e_{b}\dual e^{c}.
\]
Using \eqref{hypo2} and the properties of $(e_{j}(t))_{1\leq j\leq d}$ in Lemma \ref{lemma2.1} we obtain by a routine computation that
\beq\label{et.7}
\begin{array}{l}
\Gamma^{0}_{0b}(t)= \Gamma^{c}_{00}(t)=0,\ 
\Gamma^{0}_{ab}(t), \ \Gamma^{a}_{0b}(t)\in \cS^{-1-\mu}(\rr; \cinfb(V)),\\[2mm]
\Gamma^{c}_{ab}(t)- \Gamma^{c}_{ab}(\pm\infty)\in \cS^{-\mu}(\rr^{\pm}; \cinfb(V))\hbox{ if }a, b, c\neq 0.
\end{array}
\eeq
If we reintroduce the index $i$ and set $V= V_{i}$, then the seminorms in \eqref{et.7} are uniform with respect to $i\in \nn$.    Therefore the limits $\lim_{t\to \pm\infty}\nabla^{S}_{e_{j}(t)}$  exist and
\[
\nabla^{S}_{e_{j}(t)}- \lim_{t\to \pm\infty}\nabla^{S}_{e_{j}(t)}\in \cS^{-\mu}(\rr^{\pm}; \Psi^{1}(\Sigma; S(\Sigma))).
\]
Using \eqref{et.3} we obtain also that $ \lim_{t\to \pm\infty}\nabla^{S}_{e_{j}(t)}= \nabla^{S_{\pm\infty}}(e_{j}(\pm\infty))$, which proves {\it v)}. The proof of {\it iv)} is similar. \qed

\subsubsection{Identification by parallel transport}
  For $f\in \cinf(\Sigma_{s}; S(\Sigma_{s}))$ we denote by $\mathcal{T}(s)f= \psi$ the solution of
\beq\label{e1.8}\left\{
\begin{array}{l}
\nabla^{S}_{\p_{t}}\psi=0\hbox{ in }\rr\times \Sigma,\\
\psi_{| \Sigma_{s}}=f,
\end{array}\right.
\eeq
and set
\[
\mathcal{T}(t, s)f= \mathcal{T}(s)f\traa{ \Sigma_{t}},
\]
\beq\label{e6.2}
\begin{array}{l}
\mathcal{T}: \cinf(\rr; \cinf(\Sigma, S(\Sigma)))\to \cinf(M; S(M))\\
\psi(t)\mapsto (\mathcal{T}\psi)(t)= |{h}_{t}|^{-\frac{1}{4}}|{h}_{0}|^{\frac{1}{4}}\mathcal{T}(t, 0)\psi(t),
\end{array}
\eeq
We denote by $\nu_{0}$ the Hilbertian scalar product
\[
\bar{f}\dual \nu_{0} f\defeq \i\int_{\Sigma}\bar{f}\dual \beta_{0}\gamma_{0}(e_{0})f |{h}_{0}|^{\12}d\rx, \ f\in \coinf(\Sigma, S(\Sigma)).
\]
Using  \eqref{troup.1} we obtain the following lemma, see \cite[Lemma 6.1]{St}. 
\begin{lemma}\label{lemma3.1}
 One has
 \ben
 \item $\cT(s,t)\gamma_{t}(e_{0})\cT(t, s)= \gamma_{s}(e_{0})$, $t,s\in I$,
\item $\cT(s,t)\gamma_{t}(e_{j}(t))\cT(t, s)= \gamma_{s}(e_{j}(s))$, $t,s\in I$,
\item $\cT(t,s)^{*}\beta_{t}\cT(t,s)= \beta_{s}$, $t,s\in I$.
 \een
\end{lemma}

\subsubsection{Reduction of the Dirac operator}\label{reductionum}
\begin{proposition}\label{prop1.1}
Let
\[
\mathsf{D}\defeq \cT^{-1}{D}\cT.
\]
Then
\ben
\item the map
 \[
 \cT: (\Sol (\mathsf{D}), \nu_{0})\tosim (\Sol({D}), \nu_{0})
 \]
 is unitary.
 \item  We have 
\[
\mathsf{D}= - \gamma_{0}(e_{0})\p_{t}+\i \gamma_{0}(e_{0})\mathsf{H}(t),
\]
where $\mathsf{H}(t)\in\cS^{0}(\rr, \Psi^{1}(\Sigma, S(\Sigma)))$ has the following properties:

\medskip
 {\rm (2i)} $\sigma_{\rm pr}(\mathsf{H}(t))(\rx, k)= - \gamma_{0}(e_{0})\gamma_{0}({h}_{t}(\rx)^{-1}k)$.

\medskip
{\rm (2ii)} there exist $\mathsf{H}_{\oito}\in\Psi^{1}(\Sigma; L(S(\Sigma)))$  elliptic, formally selfadjoint for  $\nu_{0}$ such that
\[
\mathsf{H}(t)- \mathsf{H}_{\oito}\in \cS^{-\mu}(\rr^{\pm}; \Psi^{1}(\Sigma; L(S(\Sigma)))).
\]
$\mathsf{H}_{\oito}$ is selfadjoint with domain $H^{1}(\Sigma; S(\Sigma))$ and $0\not\in \sigma(H_{\pm \infty})$.

\medskip
{\rm (2iii)} $\mathsf{H}(t)$ is formally selfadjoint for $\nu_{0}$ and selfadjoint with domain $H^{1}(\Sigma; S(\Sigma))$.

\een
\end{proposition}
\proof 
(1) is obvious since $\cT(0, 0)=\one$.   We have $\cT^{-1}\gamma(e_{0})\cT= \gamma_{0}(e_{0})$ by Lemma \ref{lemma3.1} and
\begin{equation}
\label{e3.9}
\mathcal{T}\p_{t}\mathcal{T}^{-1}=\nabla_{e_{0}}^{S}- \frac{1}{4}\p_{t}|h_{t}||h_{t}|^{-1}.
\end{equation}
If we fix over some open set $U= \rr\times V$,  a local oriented and time oriented orthonormal frame $(e_{a})_{0\leq a \leq d}$  as in Lemma \ref{lemma2.1}, we have
\beq\label{e1.1bd}
\begin{array}{l}
\cT^{-1}{D}\cT= - \gamma_{0}(e_{0})\p_{t} + \i \gamma_{0}(e_{0})\mathsf{H}(t), 
\\[2mm]
\mathsf{H}(t)=  \cT^{-1}{H}(t)\mathcal{T}+ \frac{1}{4}\p_{t}|{h}_{t}||{h}_{t}|^{-1},\hbox{ for}\\[2mm]
{H}(t)\defeq  \i \gamma_{t}(e_{0})\gamma_{t}(e_{j}(t))\nabla_{e_{j}(t)}^{S}+ \i  \gamma_{t}(e_{0}){m},
\end{array}
\eeq
where in the second line we sum only over  $1\leq j \leq d$.

Let us now prove the properties of $\mathsf{H}(t)$ stated in (2).  By Prop. \ref{prop4.0} we obtain that
\begin{equation}
\label{et.9}
{H}(t)- {H}_{\oito}\in \cS^{-\mu}(\rr^{\pm}; \Psi^{1}(\Sigma; L(S(\Sigma)))),
\end{equation}
\def\oi{\outin}
for
\beq\label{troup.6}
{H}_{\oi}= \i \gamma_{\infty}(e_{0})(\gamma_{\infty}(e_{j}(\outin))\nabla_{e_{j}(\outin)}^{\oi}+ {m}_{\oi}).
\eeq
Let us now consider the maps $\cT(t,s)$. We have
\[
\nabla^{S}_{e_{0}}\psi= \p_{t}\psi+ \frac{1}{4} \Gamma^{a}_{0b} \pmb{\gamma}_{a}\pmb{\gamma}^{b}\psi= \p_{t}\psi,
\]
since $\nabla_{e_{0}}e_{a}=0$. It follows that the matrix of $\cT(t,s)$ in the bases $\cB(s)$ and $\cB(t)$ equals the identity matrix.
Using then \eqref{et.2} we obtain that the limits $\cT(\pm\infty, 0)= \lim_{t\to \pm\infty}\cT(t, 0)\in \cinfb(\Sigma; L(S(\Sigma)))$ exist and that
\begin{equation}
\label{et.10}
\cT(t, 0)- \cT(\pm\infty, 0)\in \cS^{-\mu}(\rr^{\pm}; \cinfb(\Sigma; L(S(\Sigma)))).
\end{equation}
Combining \eqref{et.9} and  \eqref{et.10}, we obtain that
\beq\label{et.11}
\mathsf{H}_{\oi}\defeq \cT(0, \pm\infty)|{h}_{\oi}|^{\frac{1}{4}}|{h}_{0}|^{-\frac{1}{4}}{H}_{\pm\infty}|{h}_{0}|^{\frac{1}{4}}|{h}_{\oi}|^{-\frac{1}{4}}\cT(\pm\infty, 0)
\eeq
belongs to $\Psi^{1}(\Sigma; L(S(\Sigma)))$ and
\begin{equation}
\label{e1.1be}
\mathsf{H}(t)- \mathsf{H}_{\oi}\in \cS^{-\mu}(\rr^{\pm}; \Psi^{1}(\Sigma; L(S(\Sigma)))).
\end{equation}
The principal symbol of $\mathsf{H}(t)$ is clearly equal to $- \gamma_{0}(e_{0})\gamma_{0}({h}_{t}(\rx)^{-1}k)$, which proves {\rm (2i)}.

Let us now prove the remaining parts of {\rm (2ii)}. From \eqref{blito}  we obtain that $\mathsf{H}_{\pm\infty}$ is the spatial part of the Dirac operator for the  static metric $g_{\pm\infty}$. Using hypothesis ${\rm (H4)}$ (and remembering that we removed the tildes) we obtain that $0\not\in \sigma(H_{\pm\infty})$. The selfadjointness  of $\mathsf{H}_{\oi}$ on $H^{1}(\Sigma; S(\Sigma))$ follows by the usual ellipticity argument. 

Finally we know that if ${D}\psi= 0$ then
\[
\int_{\Sigma}\bar{\psi}(t, \cdot)\dual \beta_{t} \gamma_{t}(e_{0}) \psi(t, \cdot)|{h}_{t}|^{\12}d\rx
\]
is independent on $t$, hence if $\tilde{\psi}= \mathcal{T}^{-1}\psi$  we obtain using
\[
\mathcal{T}\p_{t}\mathcal{T}^{-1}=\nabla_{e_{0}}^{S}- \frac{1}{4}\p_{t}|{h}_{t}||{h}_{t}|^{-1}
\]
that
\[
\int_{\Sigma}\bar{\tilde{\psi}}(t, \cdot)\dual \beta_{0} \gamma_{0}(e_{0}) \tilde{\psi}(t, \cdot)|{h}_{0}|^{\12}d\rx
\]
is independent on $t$.  Since   $\p_{t}\tilde{\psi}= \i \mathsf{H}(t)\tilde{\psi}$ this implies that $\mathsf{H}(t)= \mathsf{H}^{*}(t)$  on $\coinf(\Sigma; S(\Sigma))$ for $\nu_{0}$. The fact that (the closure of) $\mathsf{H}(t)$ is then selfadjoint on $H^{1}(\Sigma; S(\Sigma))$ follows from the standard argument, using the ellipticity of $\mathsf{H}(t)$.  \qed

\subsection{Some preparations}\label{sec4.1b}
The space $\coinf(\rr\times \Sigma; S(\Sigma))$ is equipped with the Hilbertian scalar product
\[
\bar{\psi}\dual \nu \psi= \int_{\rr\times \Sigma}\bar{\psi}\dual \beta_{0}\gamma_{0}(e_{0})\psi dt|{h}_{0}|^{\12}dx,
\]
while   $\coinf(\Sigma; S(\Sigma))$ is equipped with
\beq\label{defdeprodo}
\bar{f}\dual \nu_{0} f= \int_{\Sigma}\bar{f}\dual \beta_{0}\gamma_{0}(e_{0})f |h_{0}|^{\12}dx.
\eeq
Adjoints of operators will always be computed with respect to these scalar products. Our reference Hilbert space is
\[
\cH= L^{2}(\Sigma; S(\Sigma)),
\]
equal to the completion of $\coinf(\Sigma; S(\Sigma))$  for $\nu_{0}$.

The following lemma is the analog of \cite[Lemma 6.3]{St}.
\begin{lemma}\label{lemma4.1}
 There exists $R_{-\infty}\in \coinf(\rr; \cW^{-\infty}(\Sigma; S(\Sigma)))$ with $R_{-\infty}(t)= R_{-\infty}(t)^{*}$ and $\delta>0$ such that 
 \[
 \sigma(\mathsf{H}(t)+ R_{-\infty}(t))\cap [-\delta, \delta]= \emptyset.
 \]
 \end{lemma}
\proof we follow the proof  in \cite[Lemma 6.3]{St}. By  Prop. \ref{prop1.1} (2), we know that there exists $\delta>0$ such that $\sigma(\mathsf{H}(t))\cap [-\delta, \delta]= \emptyset$ for $|t|\gg 1$, so the modification $R_{-\infty}(t)$ can be taken compactly supported in $t$. \qed
\subsubsection{Unitary group}
Let us denote by $\mathsf{U}(t,s)$, $s,t\in I$ the unitary evolution generated by $\mathsf{H}(t)$, ie the solution of
\[
\left\{
\begin{array}{l}
\p_{t}\mathsf{U}(t, s)= \i \mathsf{H}(t)\mathsf{U}(t, s), \\
 \p_{s}\mathsf{U}(t,s)= - \i \mathsf{U}(t, s)\mathsf{H}(s),\\
\mathsf{U}s, s)= \one.
\end{array}\right.
\]
The properties of $\mathsf{H}(t)$ imply that $\mathsf{U}(t,s)$ is well-defined by a classical result of Kato, see eg \cite{SG}. 
\begin{lemma}\label{lemma5.2}
 $\mathsf{U}(t,s)$ are uniformly bounded in $B(H^{m}(\Sigma; S(\Sigma)))$ for $t,s\in \rr$, $m\in \rr$.
\end{lemma}
\proof
 Let us set
\[
S(t)\defeq (\mathsf{H}^{2}(t)+1)^{\12}, \ S_{\oi}\defeq (\mathsf{H}^{2}_{\oi}+1)^{\12}.
\]
By Prop. \ref{prop4.3} we obtain that
\beq\label{e4.10}
\begin{array}{l}
S(t)\in \cS^{0}(\rr; \Psi^{1}(\Sigma; S(\Sigma))), \\[2mm]
 S(t)- S_{\oi}\in \cS^{-\mu}(\rr^{\pm}; \Psi^{1}(\Sigma; S(\Sigma))),
\end{array}
\eeq
and $\sigma_{\rm pr}(S(t))(x, k)= (k\dual h^{-1}_{t}(\rx)k)^{\12}$. This implies that
\[
C_{m}^{-1 }\| S^{m}(t)u\|_{0}\leq \| u\|_{m}\leq C_{m}\| S^{m}(t)u\|_{0}, \ t\in \rr,
\]
 where we denote by $\| \cdot\|_{m}$ the norm in $H^{m}(\Sigma; S(\Sigma))$.

 For $u\in S^{m}(s)\coinf(\Sigma; S(\Sigma))$ we set
\[
f(t)= \|\mathsf{U}(s,t)S^{m}(t)\mathsf{U}(t, s)S^{-m}(s)u\|_{0},
\] which is finite since $\mathsf{U}(t,s)$ preserves $\coinf(\Sigma; S(\Sigma))$.  We have
\[
\begin{array}{rl}
|f'(t)|\leq& \|\mathsf{U}(s,t)\p_{t}S^{m}(t)\mathsf{U}(t,s)S^{-m}(s)u\|_{0} \\[2mm]
=&\|\mathsf{U}(s,t)\p_{t}S^{m}(t)S^{-m}(t)\mathsf{U}(t,s)\mathsf{U}(s,t)S^{m}(t)\mathsf{U}(t,s)S^{-m}(s)u\|_{0}\\[2mm]
\leq&  \|\mathsf{U}(s,t)\p_{t}S^{m}(t)S^{-m}(t)\mathsf{U}(t,s)\|_{B(\cH)}f(t)\leq Ct^{-1- \mu} f(t), 
\end{array}
\]
where we use \eqref{e4.10} in the last inequality. 
By Gronwall's inequality we obtain that $f(t)\leq C f(s)$ for $t,s\in \rr$ hence 
\[
\begin{array}{rl}
&\|S^{m}(t)\mathsf{U}(t, s)S^{-m}(s)u\|\\[2mm]
=& \|\mathsf{U}(s,t)S^{m}(t)\mathsf{U}(t, s)S^{-m}(s)u\| \leq C \| u\|,\ u\in S^{m}(s)\coinf(\Sigma; S(\Sigma)),
\end{array}
\]
 which proves the lemma since  $S^{m}(s)\coinf(\Sigma; S(\Sigma))$ is dense in $L^{2}(\Sigma; S(\Sigma))$. \qed

\subsubsection{Some preparations}\label{sec4.2.1}
We next introduce some classes of maps between pseudodifferential operators.  
These classes are similar to the ones considered in \cite[Subsect. 6.3]{St}, with the behavior for large times taken into account.

We will use the short hand notation introduced in Subsect. \ref{sec2.2} and denote $\cS^{\delta}(\rr; \Psi^{m}(\Sigma; L(S(\Sigma))))$ simply by $\cS^{\delta, m}$. We set $\cS^{\infty, \infty}= \bigcup_{\delta, m\in \rr}\cS^{\delta, m}$. 

\begin{definition}\label{def4.1}
 Let $\pmb{\delta}: \rr\to \rr$ and $p\in \rr$. We denote by $\cF_{- \pmb{\delta}, -p}$ the set of maps $F: \cS^{0, 0}\to \cS^{\infty, \infty}$ such that \[
 F:\cS^{-\mu, -1}\to \cS^{-\mu- \pmb{\delta}(\mu), -p}, \ \forall \mu>0,
 \]
and:
\[
R_{1}- R_{2}\in \cS^{-\mu- \epsilon, -1-j}\Rightarrow  F(R_{1})- F(R_{2})\in \cS^{- \mu -\pmb{\delta}(\mu)-\epsilon, -p- j}, \ \forall \epsilon>0, j\in \nn.
\]
\end{definition}
An element of $\cF_{-\pmb{\delta}, -p}$ will be denoted by $F_{-\pmb{\delta}, -p}$.  The following proposition is the analog of \cite[Lemma A.1]{GW1}, \cite[Prop. 6.6]{St}. It is an abstract formulation of an  ubiquitous argument in pseudodifferential calculus, consisting in solving  recursive equations to determine successive terms  in the  symbolic expansion of a pseudodifferential operators. \begin{proposition}\label{prop5.1}
 Let $A\in\cS^{- \mu_{1}, -1}$, $\mu_{1}>0$ and $F_{0, -2}\in \cF_{0, -2}$. Then there exists a solution $R\in \cS^{- \mu_{1}, -1}$, unique modulo $\cS^{- \mu_{1}, -\infty}$ of the equation:
 \[
 R= A+ F_{0, -2}(R)\hbox{ mod }\cS^{- \mu_{1}, -\infty}.
 \]
 \end{proposition}
 \proof Let us denote $F_{0, -2}$ simply by $F$. We set $S_{0}= A$, $S_{n}= A+ F(S_{n-1})$ for $n\geq 1$. We have $S_{1}- S_{0}= F(A)$ and $S_{n}- S_{n-1}= F(S_{n-1})- F(S_{n-2})$. Since $F\in \cF_{0, -2}$ we obtain by induction that $S_{n}- S_{n-1}\in \cS^{- \mu_{1}, -(n+1)}$. We take $R\in \cS^{-\mu_{1}, -1}$ such that $R\sim S_{0}+ \sum_{0}^{\infty}S_{n}- S_{n-1}$ which solves the equation modulo $\cS^{- \mu_{1}, -\infty}$. If $R_{1}, R_{2}$ are two solutions then
 $R_{1}- F_{2}= F(F_{1})- F(R_{2})$ modulo $\cS^{-\mu_{1}, -\infty}$ hence using that $F\in \cF_{0, -2}$ we obtain by induction on $n$ that $R_{1}- R_{2}\in \cS^{-\mu_{1}, -n}$ for all $n\in \nn$ which proves uniqueness modulo $\cS^{-\mu_{1}, -\infty}$. \qed

We now collect some useful properties of the sets $\cF_{-\pmb{\delta},-p}$.
\begin{lemma}\label{lemma4.2}
  \ben
 \item If $A\in\cS^{-\varrho, k}$  and $F_{-\pmb{\delta}, -p}\in \cF_{-\pmb{\delta},-p}$, then  the maps \[
 \begin{array}{l}
 A F_{-\pmb{\delta}, -p}: R\mapsto AF_{-\pmb{\delta}, -p}(R),\\[2mm]
 F_{-\pmb{\delta},-p}A: R\mapsto F_{-\pmb{\delta},-p}(R)A
 \end{array}
 \]
  belong to $\cF_{-\pmb{\delta}-\varrho,-p+k}$ for $k\leq p$.
 \item If $F_{-\pmb{\delta}_{i},-p_{i}}\in \cF_{-\pmb{\delta}_{i}, -p_{i}}$ then  the map 
 \[
 F_{-\pmb{\delta_{1}}, -p_{1}}F_{-\pmb{\delta}_{2}, -p_{2}}: R\mapsto F_{-\pmb{\delta_{1}}, -p_{1}}(R)F_{-\pmb{\delta_{2}}-p_{2}} (R)
 \] 
 belongs to $\cF_{-\pmb{\delta_{1}}- \pmb{\delta_{2}}-\pmb{\mu},-p_{1}- p_{2}}$, where $\pmb{\mu}$ is the map $\mu\mapsto \mu$.
  \item the map $R\to R^{p}$ belongs to $\cF_{-(p-1)\pmb{\mu}, -p}$ for $p\in \nn^{*}$.
 \item the map $R\mapsto \e^{R}$ belongs to $\cF_{0,0}$.
\item one has $\e^{R}= 1 + R + F_{-\pmb{\mu}, -2}(R)$, where $F_{-\pmb{\mu}, -2}\in \cF_{-\pmb{\mu}, -2}$.
 \een

\end{lemma}
\proof (1) and (2) are easy.  We check that $R\mapsto R$ belongs to $\cF_{0, -1}$ and use then (2) to obtain (3).  To prove (4) we write $\e^{R}= \sum_{n\geq 0}\frac{1}{n!}R^{n}$ and obtain that $\e^{R}\in \cS^{0, 0}$ if $R\in \cS^{-\mu, -1}$.  We have
\[
\e^{R_{1}}- \e^{R_{2}}= \int_{0}^{1}\e^{ \theta R_{1}}(R_{1}- R_{2})\e^{(1- \theta)R_{2}}d\theta
\]
and obtain that $\e^{R_{1}}- \e^{R_{2}}\in \cS^{- \mu- \epsilon, -1- j}$ if $R_{1}- R_{2}\in \cS^{-\mu-\epsilon, -1-j}$, which completes the proof of (4). To prove (5) we write

\[
\e^{R}= \one + R+\int_{0}^{1}(1- \theta)R^{2}\e^{ \theta R}d\theta\eqdef 1+ R+ F(R)
\]
and obtain by (2), (3) and (4) that $F\in \cF_{\pmb{\mu}, -2}$. \qed

\subsection{Construction of some projections}\label{sec4.2}
We now follow the constructions in \cite[Subsect. 6.4]{St}, adapting the results to our  framework.

\begin{proposition}\label{prop4.1}
 There exist time-dependent projections \[
 P^{\pm}(t)\in \cS^{0}(\rr; \Psi^{0}(\Sigma; L(S(\Sigma))))
 \]
   and time-dependent operators \[
   R(t)\in\cS^{-1-\mu}(\rr; \Psi^{-1}(\Sigma; L(S(\Sigma))))
   \] such that
 \ben
 \item $P^{\pm}(t)= P^{\pm}(t)^{*}$, $P^{+}(t)+ P^{-}(t)=\one$;
 \item $P^{\pm}(t)- \one_{\rr^{\pm}}(\mathsf{H}_{\oi})\in \cS^{-\mu}(\rr^{\pm}, \Psi^{0}(\Sigma; L(S(\Sigma))))$;
 \item $R(t)= R^{*}(t)$;
 \item $\p_{t}P^{\pm}(t)+ [P^{\pm}(t),\i \tilde{\mathsf{H}}(t)]\in \cS^{-1-\mu}(\rr; \Psi^{-\infty}(\Sigma; L(S(\Sigma))))$ for
 \[
\tilde{\mathsf{H}}(t)= \e^{\i R(t)}\mathsf{H}(t)\e^{- \i R(t)}+ \i^{-1}\p_{t}\e^{\i R(t)}\e^{- \i R(t)};
 \]
 \item 
 \[
 WF( \mathsf{U}(\cdot,0)\e^{-\i R(0)}P^{\pm}(0)\e^{\i R(0)})'\subset (\cN^{\pm}\cup \cF)\times T^{*}\Sigma,
 \]
 for $\cF= \{k=0\}\subset T^{*}M$.
 \een
\end{proposition}
\proof we follow the proof of \cite[Prop.6.8]{St}, taking into account the time decay of the various operators. 

{\it Step1.} In Step 1 we replace $\mathsf{H}(t)$ by $\hat{\mathsf{H}}(t)= \mathsf{H}(t)+R_{-\infty}(t)$ as in Lemma \ref{lemma4.1}. Let $\hat{\mathsf{U}}(t,s)$ the unitary group with generator $\hat{\mathsf{H}}(t)$. From Lemma \ref{lemma5.2} and Duhamel's formula we obtain that $\mathsf{U}(t,s)- \hat{\mathsf{U}}(t,s)\in \cinfb(\rr^{2}; \cW^{-\infty}(\Sigma; L(S(\Sigma))))$ so we can replace $\mathsf{H}(t)$ by $\hat{\mathsf{H}}(t)$. Denoting $\hat{\mathsf{H}}(t)$ again by $\mathsf{H}(t)$ we can assume without loss of generality that $[-\delta, \delta]\cap \sigma(\mathsf{H}(t))= \emptyset$ for $t\in \rr$.

By Prop. \ref{prop4.3} the projections
\[
P^{\pm}(t)= \one_{\rr^{\pm}}(\mathsf{H}(t))
\]
are well defined, selfadjoint with 
\beq\label{e4.11}
\begin{array}{l}
P^{\pm}(t)\in \cS^{0}(\rr; \Psi^{0}(\Sigma; L(S(\Sigma)))), \\[2mm]
 P^{\pm}(t)- \one_{\rr^{\pm}}(\mathsf{H}_{\oi})\in \cS^{-\mu}(\rr^{\pm}; \Psi^{0}(\Sigma; L(S(\Sigma)))),
\end{array}
\eeq 
so properties (1) and (2) are satisfied. We have also
\begin{equation}
\label{e5.0}
\sigma_{\rm pr}(P^{\pm})(t, \rx, k)= \one_{\rr^{\pm}}(\sigma_{\rm pr}(H)(t, \rx, k)).
\end{equation}
Since $\sigma_{\rm pr}(\mathsf{H}(t, \rx, k))= -\gamma_{0}\gamma({h}_{t}^{-1}(\rx)k)$, we obtain using the Clifford relations that:
\[
\sigma_{\rm pr}(P^{\pm})(t, \rx, k)\sigma_{\rm pr}(\mathsf{H})(t, \rx, k)= \pm\epsilon(t, \rx, k)\sigma_{\rm pr}(P^{\pm})(t, \rx, k),
\]
for $\epsilon(t, x, k)= (k\dual {h}_{t}^{-1}(\rx)k)^{\12}$. By symbolic calculus this implies that
\begin{equation}
\label{e5.-0}
P^{\pm}(t)\mathsf{H}(t)= \pm \epsilon(t, \rx, D_{x}) P^{\pm}(t)+ R_{0}^{\pm}(t), 
\end{equation}
where $R_{0}^{\pm}(t)\in \cS^{0}(\rr; \Psi^{0}(\Sigma; S(\Sigma)))$.

 {\it Step 2.} In Step 2 we  find $R(t)$ such that (4) is satisfied.
For ease of notation we denote simply by $A$ a time-dependent pseudodifferential operator $A(t)$.  By Lemma \ref{lemma4.2} we obtain easily that for $\tilde{\mathsf{H}}(t)$ defined in (4) we have:
\beq\label{e5.3b}
\tilde{\mathsf{H}}= \mathsf{H} + [R, \i \mathsf{H}]+ F_{-\pmb{\mu}, -1}(R).
\eeq
We will look for $R$ of the form
\beq
\label{e5.4}
R=T(S)=  P^{+}SP^{+}+ P^{-}S^{*}P^{-}, \ S\in\cS^{0, -1}.
\end{equation}
Note that if $F_{-\pmb{\delta}, -p}\in \cF_{-\pmb{\delta}, -p}$ then the map $S\mapsto F_{-\pmb{\delta}, -p}(T(S))$ belongs also to $\cF_{-\pmb{\delta}, -p}$ (note that the map $S\mapsto S^{*}$ belongs to  $\cF_{0, -1}$).

Since $P^{\pm}$ are projections we have
 \[
 \begin{array}{rl}
& \p_{t}P^{\pm}+ [P^{\pm}, \i \tilde{\mathsf{H}}]\\[2mm]
=& P^{+}(\p_{t}P^{\pm}+ [P^{\pm},  \i\tilde{\mathsf{H}}])P^{-}+ P^{-}(\p_{t}P^{\pm}+ [P^{\pm}, \i  \tilde{\mathsf{H}}])P^{+}.
\end{array}
 \]
 Since the second term in the rhs above is the adjoint of the first, it suffices to find $S$ such that 
 \begin{equation}
 \label{e5.5}
 P^{+}(\p_{t}P^{+}+ [P^{+}, \i \tilde{\mathsf{H}}])P^{-}\in\cS^{-1-\mu, -\infty}.
 \end{equation}
 Using \eqref{e5.3b} we obtain since $[P^{\pm}, \mathsf{H}]=0$:
 \[
 \begin{array}{rl}
 &P^{+}\left(\p_{t}P^{+}+ [P^{+}, \i \tilde{\mathsf{H}}]\right)P^{-}\\[2mm]
 =&P^{+}\left(\p_{t}P^{+}+ P^{+}\mathsf{H}P^{+}S- S P^{-}\mathsf{H} P^{-}+ F_{-\pmb{\mu}, -1}(S)\right)P^{-}
 \end{array}
 \]
 We use now \eqref{e5.-0} denoting  the scalar operator $\epsilon(t, \rx, D_{x})+ m^{2}$ for $m\gg 1$ simply by $\epsilon$ and obtain:
 \[
 \begin{array}{rl}
& P^{+}\mathsf{H}P^{+}S- S P^{-}\mathsf{H} P^{-}= \epsilon S+ S \epsilon+ R_{0}^{+}S- S R_{0}^{-}\\[2mm]
=&2 \epsilon S+ [S, \epsilon]+ R_{0}^{+}S- S R_{0}^{-}.
 \end{array}
 \]
 The maps $S\mapsto R_{0}^{+}S$, $S\mapsto SR_{0}^{-}$ belong to $\cF_{0,-1}$ by Lemma \ref{lemma4.2}, as the map $S\mapsto [\epsilon, S]$, since $\epsilon$ is scalar. 

Therefore the equation \eqref{e5.5} can be rewritten as
\[
\p_{t}P^{+}+ 2\epsilon S+ F_{0, -1}(S)\in \cS^{-1-\mu, -\infty},
\]
or equivalently as
\begin{equation}
\label{e5.6}
S+(2\epsilon)^{-1}\p_{t}P^{+}+-F_{0, -2}(S)\in  \cS^{-1-\mu, -\infty}.
\end{equation}
where $F_{0, -2}: S\mapsto - (2 \epsilon)^{-1}F_{0, -1}(S)$ belongs to $\cF_{0, -2}$. We apply Prop. \ref{prop5.1} to solve \eqref{e5.6}. 
We note that $- (2\epsilon)^{-1}\p_{t}P^{+}\in \cS^{-1-\mu, -1}$ and we find $S\in \cS^{-1-\mu, -1}$ such that
\[
\p_{t}P^{+}+ 2 \epsilon S+ F_{0, -1}(S)\in \cS^{-1-\mu, -\infty}
\]
and hence
\[
\p_{t}P^{+}+ [P^{+}, \i \tilde{\mathsf{H}}]= R_{-\infty}\in \cS^{-1-\mu, -\infty}.
\]
We have hence proved (4). Finally (5) is proved exactly as in \cite[Prop. 6.8]{St}. \qed
 \subsection{The in/out vacua for $D$}\label{sec4.3}
 In this subsection we denote again by $\tilde{g}$, $\tilde{D}$ etc the objects obtained from $g$, $D$ by conformal transformation. For example the scalar product in \eqref{defdeprodo} is now denoted $\tilde{\nu}_{0}$.

  We recall also from \ref{cauchy-evol} that $U(t, s)$ is the Cauchy evolution for $D$ associated to the foliation $(\Sigma_{t})_{t\in \rr}$.  We denote by $L^{2}(\Sigma_{t}; S(\Sigma))$ the completion of $\coinf(\Sigma_{t}; S(\Sigma))$ for the scalar product
 \[
 \bar{f}\dual \nu_{t}f= \i \int_{\Sigma_{t}} \bar{f}\dual \beta \gamma(n) dVol_{h_{t}}.
 \]
 By the facts recalled in \eqref{pre-hilbert} $U(t,s): L^{2}(\Sigma_{s}; S(\Sigma))\to L^{2}(\Sigma_{t}; S(\Sigma))$ is unitary.
 We denote by $\blibli{U}(t, s)$ the analogous Cauchy evolution for $\blibli{D}$.  By \eqref{troup.5} we have:
 \begin{equation}
 \label{troup.16}
 U(t,s)= c^{\frac{1-n}{2}}\blibli{U}(t, s)c^{\frac{n-1}{2}}.
 \end{equation}

 From the definition \eqref{troup.6} of $\blibli{H}_{\pm\infty}$ we obtain that
 the asymptotic Dirac operator $\blibli{D}_{\oi}$ associated to the ultrastatic metric $\blibli{g}_{\oi}$ equals
\[
\blibli{D}_{\oi}= - \tilde{\gamma}_{\oi}(\tilde{e}_{0})(\p_{t}-\i\blibli{H}_{\oi}).
\]

Recalling that $D_{\oi}$ are the asymptotic Dirac operators associated to $g_{\oi}$ we have as in Subsect. \ref{secvac}:
\[
D_{\oi}= -c_{\oi}^{-1} \gamma_{\pm\infty}(e_{0})(\p_{t}- \i H_{\oi}),
\]
and 
\begin{equation}
\label{troup.17}
H_{\oi}= c_{\oi}^{\frac{1-n}{2}}\blibli{H}_{\oi}c_{\oi}^{\frac{n-1}{2}}.
\end{equation}
We first consider the operator $\mathsf{D}$ in \ref{reductionum}. 
\begin{proposition}\label{propot}
Assume hypotheses ${\rm (Hi)}$, $1\leq i\leq 4$.  Then:
\ben
\item  the norm limits:
\[
\mathsf{P}^{\pm}_{{\rm out/in}}= \lim_{t\to \pm\infty}\mathsf{U}(0, t)\one_{\rr^{\pm}}(\mathsf{H}_{\oi})\mathsf{U}(t, 0)\hbox{ exist}.
\]
\item  $\mathsf{P}^{\pm}_{{\rm out/in}}$ are selfadjoint projections for $\tilde{\nu}_{0}$ with  $\mathsf{P}^{+}_{{\rm out/in}}+ \mathsf{P}^{-}_{{\rm out/in}}=\one$.
\item \[
\WF( \mathsf{U}(\cdot,0)\mathsf{P}^{\pm}_{{\rm out/in}})'\subset (\cN^{\pm}\cup \cF)\times T^{*}\Sigma
\]
for $\cF= \{k=0\}\subset T^{*}M$.
\een
\end{proposition}

\proof 
Let $P^{\pm}(t), R(t)$ be the operators constructed in Prop. \ref{prop4.1}. Setting
\[
\tilde{P}^{\pm}(t)\defeq \e^{-\i R(t)}P^{\pm}(t)\e^{\i R(t)}, \ \tilde{\mathsf{U}}(t, s)\defeq \e^{\i R(t)}\mathsf{U}(t,s)\e^{- \i R(s)},
\]
we see that $\tilde{\mathsf{U}}(t,s)$ is a strongly continuous unitary group with generator
\[
\tilde{\mathsf{H}}(t)= \e^{\i R(t)}\mathsf{H}(t)\e^{- \i R(t)}+ \i^{-1}\p_{t}\e^{\i R(t)}\e^{- \i R(t)}.
\]
Since $P^{\pm}(t)- \one_{\rr^{\pm}}(\mathsf{H}_{\oi})$ and $R(t)$ are $O(t^{-\mu})$ 
 in norm we have $\tilde{P}^{\pm}(t)- \one_{\rr^{\pm}}(\mathsf{H}_{\oi})\in O(t^{-\mu})$ and hence 
 \beq\label{troup.14}
 \mathsf{U}(0, t)\one_{\rr^{\pm}}(\mathsf{H}_{\oi})\mathsf{U}(t, 0)= \mathsf{U}(0, t)\tilde{P}^{\pm}(t)\mathsf{U}(t, 0)+ O(t^{-\mu}).
 \eeq
Next
\[
\begin{array}{rl}
&\mathsf{U}(0,t)\tilde{P}^{\pm}(t)\mathsf{U}(t, 0)=\e^{- \i R(0)}\tilde{\mathsf{U}}(0, t)\e^{\i R(t)}\tilde{P}^{\pm}(t)\e^{- \i R(t)}\tilde{\mathsf{U}}(t, 0)\e^{\i R(0)}\\[2mm]
=&\e^{- \i R(0)}\tilde{\mathsf{U}}(0, t)P^{\pm}(t)\tilde{\mathsf{U}}(t, 0)\e^{\i R(0)},
\end{array}
\]
and
\beq\label{troup.15}
\begin{array}{rl}
&\p_{t}\left(\tilde{\mathsf{U}}(0, t)P^{\pm}(t)\tilde{\mathsf{U}}(t, 0)\right)\\[2mm]
=& \tilde{\mathsf{U}}(0, t)\left(\p_{t}P^{\pm}(t)+ [P^{\pm}(t), \i \tilde{\mathsf{H}}(t)]\right)\tilde{\mathsf{U}}(t, 0)= \tilde{\mathsf{U}}(0, t)R_{-\infty}(t)\tilde{\mathsf{U}}(t, 0),
\end{array}
\eeq
where $R_{-\infty}(t)\in \cS^{-1-\mu, -\infty}$, by Prop. \ref{prop4.1}. Therefore  by \eqref{troup.14}, \eqref{troup.15} the limit
\[
\mathsf{P}^{\pm}_{{\rm out/in}}= \lim_{t\to \pm\infty}\mathsf{U}(0, t)\one_{\rr^{\pm}}(\mathsf{H}_{\oi})\mathsf{U}(t, 0) \hbox{ exists }
\]
and
\[
\mathsf{P}^{\pm}_{{\rm out/in}}= \tilde{P}^{\pm}(0)\pm\int_{0}^{\pm\infty} \tilde{\mathsf{U}}(0, t)R_{-\infty}(t)\tilde{\mathsf{U}}(t, 0)dt= \tilde{P}^{\pm}(0)+ R_{\pm\infty},
\]
where $R_{\pm\infty}\in \Psi^{-\infty}$, using the uniform estimates in Lemma \ref{lemma5.2} .   $\mathsf{P}^{\pm}_{{\rm out/in}}$ are clearly selfadjoint projections for $\tilde{\nu}_{0}$ with $\mathsf{P}^{+}_{{\rm out/in}}(0)+ \mathsf{P}^{-}_{{\rm out/in}}(0)= \one$.

From  Prop. \ref{prop4.1} (5) and the fact that 
$\mathsf{P}^{\pm}_{{\rm out/in}}- \tilde{P}^{+}(0)$ is a smoothing operator we obtain  that $\WF(\mathsf{U}(\cdot, 0)\mathsf{P}^{\pm}_{{\rm out/in}})'\subset (\cN^{\pm}\cup \cF)\times T^{*}\Sigma$.  \qed

Next we consider the Dirac operator $\tilde{D}$.
\begin{proposition}\label{prop-main}
Assume hypotheses ${\rm (Hi)}$, $1\leq i\leq 4$.   Then 
 \ben
 \item the norm limits
 \[
 \blibli{P}^{\pm}_{\rm out/in}= \lim_{t\to \pm\infty}\blibli{U}(0, t)\one_{\rr^{\pm}}(\blibli{H}_{\pm\infty})\blibli{U}(t,0) \hbox{ exist}.
 \]
 \item $\blibli{P}^{\pm}_{\oi}$ are selfadjoint projections for  the scalar product $\tilde{\nu}_{0}$ with $\blibli{P}^{+}_{\oi}+ \blibli{P}^{+}_{\oi}= \one$.
 \item \[
 \WF(\blibli{U}(\cdot, 0)\blibli{P}^{\pm}_{\oi})'\subset (\cN^{\pm}\cup \cF)\times T^{*}\Sigma.
 \]
 \een
\end{proposition}
\proof
 We obtain easily from \eqref{e6.2}
that
\beq\label{troup.8}
\mathsf{U}(t,s)= \cT(0, t)|\blibli{h}|_{0}^{-\frac{1}{4}}|\blibli{h}|_{t}^{\frac{1}{4}} \blibli{U}(t,s)|\blibli{h}|_{0}^{\frac{1}{4}}|\blibli{h}|_{s}^{-\frac{1}{4}}\cT(s, 0).
\eeq
 This implies that
\[
\begin{array}{rl}
&\mathsf{U}(0, t)\one_{\rr^{\pm}}(\mathsf{H}_{\oi})\mathsf{U}(t, 0)\\[2mm]
=&\blibli{U}(0, t)|\blibli{h}|_{0}^{\frac{1}{4}}|\blibli{h}|_{t}^{-\frac{1}{4}}\cT(t, 0)\one_{\rr^{\pm}}(\mathsf{H}_{\oi})\cT(0, t)|\blibli{h}|_{0}^{-\frac{1}{4}}|\blibli{h}|_{t}^{\frac{1}{4}}\blibli{U}(t, 0).
\end{array}
\]
By \eqref{et.10}  $\cT(0, t)- \cT(0, \pm\infty)\in O(t^{-\mu})$ in norm and $|h|_{t}- |h|_{\oi}\in O(t^{-\mu})$, hence
\[
\begin{array}{rl}
&\mathsf{U}(0, t)\one_{\rr^{\pm}}(\mathsf{H}_{\oi})\mathsf{U}(t, 0)\\[2mm]
=& \blibli{U}(0, t)|\blibli{h}|_{0}^{\frac{1}{4}}|\blibli{h}|_{\oi}^{-\frac{1}{4}}\cT(\pm\infty, 0)\one_{\rr^{\pm}}(H_{\oi})\cT(0, \pm\infty)|\blibli{h}|_{0}^{-\frac{1}{4}}|\blibli{h}|_{\oi}^{\frac{1}{4}}\blibli{U}(t, 0)+ O(t^{-\mu})\\[2mm]
=& \blibli{U}(0, t)\one_{\rr^{\pm}}(\blibli{H}_{\oi})\blibli{U}(t, 0)+ O(t^{-\mu}).
\end{array}
\]
where in the last line we use \eqref{et.11} and the fact that $\cT(0, \pm\infty)|\blibli{h}|_{0}^{-\frac{1}{4}}|\blibli{h}|_{\oi}^{\frac{1}{4}}$ is unitary for the scalar products $\tilde{\nu}_{\oi}$ and $\tilde{\nu}_{0}$.  

Therefore  the norm limit in (1) exist and equal the projections $\mathsf{P}^{\pm}_{\oi}$ in Prop. \ref{propot}.   This also implies (2). (3) follows from \eqref{troup.8} and the analogous statement in Prop. \ref{propot}. \qed

Finally we prove the main result of this paper.
 \begin{theoreme}\label{mainmaino}
Assume hypotheses ${\rm (Hi)}$, $1\leq i\leq 4$.  Then:
\ben
\item the norm limits 
\[
P^{\pm}_{\oi}= \lim_{t\to \pm\infty}U(0, t)\one_{\rr^{\pm}}(H_{\oi})U(t, 0)\hbox{ exist}.
\]
\item if 
\[
\lambda^{\pm}_{\oibis}=  \i \gamma(n)P^{\pm}_{\oi}
\]
$\lambda^{\pm}_{\oibis}$ are the Cauchy surface covariances of a pure Hadamard state  for $D$ $\omega_{\oibis}$ called the {\em out/in vacuum state}.
\een
\end{theoreme}
\proof 
Let us denote by $c_{t}$ the restriction of the conformal factor $c$ to $\Sigma_{t}$.
From \eqref{troup.5} we obtain that 
\[
U(t,s)= c_{t}^{\frac{1-n}{2}}\tilde{U}(t,s)c_{s}^{\frac{n-1}{2}}, \ t,s\in \rr,
\]
and hence
\[
\begin{array}{rl}
&U(0, t)\one_{\rr^{\pm}}(H_{\oi})U(t, 0)\\[2mm]
=& c_{0}^{\frac{1-n}{2}}\tilde{U}(0, t)c_{t}^{\frac{n-1}{2}}\one_{\rr^{\pm}}(H_{\oi})c_{t}^{\frac{1-n}{2}}\tilde{U}(t, 0)c_{0}^{\frac{n-1}{2}}\\[2mm]
=& c_{0}^{\frac{1-n}{2}}\tilde{U}(0, t)c_{\oi}^{\frac{n-1}{2}}\one_{\rr^{\pm}}(H_{\oi})c_{\oi}^{\frac{1-n}{2}}\tilde{U}(t, 0)c_{0}^{\frac{n-1}{2}}+ O(t^{-\mu})\\[2mm]
=& c_{0}^{\frac{1-n}{2}}\tilde{U}(0, t)\one_{\rr^{\pm}}(\tilde{H}_{\oi})\tilde{U}(t, 0)c_{0}^{\frac{n-1}{2}}+ O(t^{-\mu})
\end{array}
\]
since $H_{\oi}= c_{\oi}^{\frac{1-n}{2}}\tilde{H}_{\oi}c_{\oi}^{\frac{n-1}{2}}$. By Prop. \ref{prop-main} we obtain that 
\beq\label{etroup.100}
P^{\pm}_{\oi}= \lim_{t\to \pm\infty}U(0, t)\one_{\rr^{\pm}}(H_{\oi})U(t, 0)= c_{0}^{\frac{1-n}{2}}\tilde{P}^{\pm}_{\oi}c_{0}^{\frac{n-1}{2}}.
\eeq
By Prop. \ref{prop20.1}  we obtain that $(c_{0}^{\frac{n-1}{2}})^{*}\tilde{\nu}_{0}c_{0}^{\frac{n-1}{2}}= \nu_{0}$.  By Prop. \ref{prop-main} $P^{\pm}_{\oi}$ are hence selfadjoint projections for $\nu_{0}$ with $P^{+}_{\oi}+ P^{-}_{\oi}= \one$. Therefore $\lambda^{\pm}_{\oibis}$ are the Cauchy surface covariances of  pure Hadamard states $\omega_{\oibis}$ for $D$. 

Finally \eqref{etroup.100} and Prop. \ref{prop-main} (3) imply that $\WF(U(\cdot, 0)P^{\pm}_{\oibis})'\subset (\cN^{\pm}\cup \cF)\times T^{*}\Sigma$. Since $\cF\cap \cN= \emptyset$ we obtain by Prop. \ref{prop15.0a} that $\omega_{\oibis}$ are Hadamard states. \qed

\end{document}